\documentclass[12pt]{amsart}

\usepackage{amssymb}
\usepackage{ifthen}
\usepackage{mathrsfs}
\usepackage{pgfpages}
\usepackage{pict2e}
\usepackage{xargs}
\usepackage{xspace}

\newcommand{\definedsymbol}[1]{$#1$}
\newcommand{\definedterm}[1]{\emph{#1}}

\newcommand{\Bairespace}[1][]{\ifthenelse{\equal{#1}{}}{\functions{\N}{\N}}{\functions{#1}{\N}}}
\newcommand{\Bairetree}[1][]{\ifthenelse{\equal{#1}{}}{\functions{<\N}{\N}}{\functions{<#1}{\N}}}
\newcommand{\branches}[1]{[#1]}

\newcommand{\calE}{\mathscr{E}}
\newcommand{\calF}{\mathscr{F}}

\newcommand{\calN}{\mathcal{N}}

\newcommand{\calU}{\mathcal{U}}

\newcommandx{\CantorCantorspace}[2][1 =, 2 =]{
  \ifthenelse{\equal{#1}{}}{\functions{\N}{(\Cantorspace)}}{
    \ifthenelse{\equal{#2}{}}{\functions{#1}{(\Cantorspace)}}{\functions{#2}{(\Cantorspace[#1])}}
  }
}
\newcommandx{\CantorCantortree}[2][1 =, 2 =]{
  \ifthenelse{\equal{#1}{}}{\functions{<\N}{(\Cantortree)}}{\functions{<#1}{(\Cantortree[#2])}}
}
\newcommand{\Cantorspace}[1][]{\ifthenelse{\equal{#1}{}}{\functions{\N}{2}}{\functions{#1}{2}}}
\newcommand{\Cantortree}[1][]{\ifthenelse{\equal{#1}{}}{\functions{<\N}{2}}{\functions{<#1}{2}}}

\newcommand{\chromaticnumber}[2][]{\chi_B(#2)}

\newcommand{\composition}{\circ}
\newcommandx{\concatenation}{\smallfrown}

\newcommand{\continuousfunctions}[2]{C(#1, #2)}

\newcommand{\domain}[1]{\mathrm{dom}(#1)}

\newcommand{\equivalenceclass}[2]{[#1]_{#2}}
\newcommand{\extendedby}{\sqsubseteq}
\newcommand{\extensions}[1]{\calN_{#1}}
\newcommand{\Eone}{\mathbb{E}_1}
\newcommand{\Etail}[1][]{\ifthenelse{\equal{#1}{}}{\mathbb{E}_t}{\mathbb{E}_t(#1)}}
\newcommand{\Ezero}[1][]{\ifthenelse{\equal{#1}{}}{\mathbb{E}_0}{\mathbb{E}_0(#1)}}
\newcommandx{\F}[3][1 =, 2 =, 3 =]{
  \ifthenelse{\equal{#1}{}}{\mathbb{F}}{
    \ifthenelse{\equal{#2}{}}{\mathbb{F}_{#1}}{
      \ifthenelse{\equal{#3}{}}{\mathbb{F}_{#1}(\CantorCantorspace[#2])}{
        \mathbb{F}_{#1}(\CantorCantorspace[#2][#3])
      }
    }
  }
}

\newcommand{\forcomeagerlymany}{\forall^*}
\newcommand{\Fsigma}{F_\sigma}
\newcommand{\from}{\colon}
\newcommand{\functiongraph}[1]{G_{#1}}
\newcommand{\functions}[2]{#2^{#1}}
\newcommandx{\G}[2][1 =, 2 =]{\ifthenelse{\equal{#2}{}}{\mathbb{G}_{#1}}{\mathbb{G}_{#1}^{#2}}}
\newcommand{\Gdelta}{G_\delta}

\newcommand{\Gzero}{\mathbb{G}_0}
\renewcommandx{\H}[2][1 =, 2 =]{\ifthenelse{\equal{#2}{}}{\mathbb{H}_{#1}}{\mathbb{H}_{#1}^{#2}}}

\newcommand{\heightcorrection}[1]{\raisebox{0pt}[0pt][0pt]{#1}}
\newcommand{\horizontalconcatenation}{\concatenation}
\newcommand{\horizontalsection}[2]{#1^{#2}}

\renewcommand{\iff}{\Longleftrightarrow}
\newcommand{\image}[2]{#1(#2)}
\renewcommand{\impliedby}{\Longleftarrow}
\renewcommand{\implies}{\Longrightarrow}

\newcommandx{\insertion}[2][2=]{\ifthenelse{\equal{#2}{}}{i_{#1}}{i_{#1}(#2)}}

\newcommandx{\intersection}[2][1 = undefined, 2 = undefined]{
  \ifthenelse{\equal{#1}{undefined}}{\cap}{
    \ifthenelse{\equal{#2}{undefined}}{\bigcap #1}{\bigcap_{#1} #2}
  }
}
\newcommandx{\interval}[3][3 =]{[#1, #2]_{#3}}
\newcommand{\Ksigma}{K_\sigma}

\newcommand{\mathand}{\text{ and }}

\newcommand{\mathor}{\text{ or }}
\newcommand{\N}{\mathbb{N}}

\newcommand{\pair}[2]{(#1, #2)}
\newcommandx{\Piclass}[2][1 =, 2 =]{
  \ifthenelse{\equal{#1}{}}{\mathbf{\Pi}}{
    \ifthenelse{\equal{#2}{}}{\mathbf{\Pi}_{#1}}{\mathbf{\Pi}^{#1}_{#2}}
  }
}

\newcommand{\preimage}[2]{#1^{-1}(#2)}
\newcommand{\projection}[1]{\mathrm{proj}_{#1}}
\newcommand{\pushforward}[2]{#1_* #2}
\newcommand{\Q}{\mathbb{Q}}
\newcommand{\quadruple}[4]{\sequence{#1, #2, #3, #4}}

\newcommandx{\R}[2][1=,2=]{\ifthenelse{\equal{#1}{}}{\mathbb{R}}{\mathbb{R}_{#1}(#2)}}
\newcommand{\relationpower}[2]{#1^{(#2)}}
\renewcommand{\restriction}[2]{#1 \upharpoonright #2}
\newcommand{\saturation}[2]{[#1]_{#2}}
\newcommandx{\sequence}[2][2 = undefined]{
  \ifthenelse{\equal{#2}{undefined}}{(#1)}{
    (#1)_{#2}
  }
}

\newcommandx{\set}[2][2 = undefined]{
  \ifthenelse{\equal{#2}{undefined}}{\{ #1 \}}{
    \{ #1 \suchthat #2 \}
  }
}
\newcommand{\setcomplement}[1]{\mathord{\sim} #1}
\newcommandx{\sets}[4][3 = undefined, 4 = undefined]{
  \ifthenelse{\equal{#3}{undefined}}{[#2]^{#1}}{
    \ifthenelse{\equal{#4}{undefined}}{[#2]^{#1}_{#3}}{[#2]^{#1}_{#3 / #4}}
  }
}
\newcommand{\sigmaclass}[1]{\sigma(#1)}
\newcommandx{\Sigmaclass}[2][1 =, 2 =]{
  \ifthenelse{\equal{#1}{}}{\mathbf{\Sigma}}{
    \ifthenelse{\equal{#2}{}}{\mathbf{\Sigma}_{#1}}{\mathbf{\Sigma}^{#1}_{#2}}
  }
}
\newcommand{\singleton}[1]{\set{#1}}

\newcommand{\suchthat}{\mid}
\newcommand{\singletonsequence}[1]{\sequence{#1}}

\newcommand{\terminal}[1]{T^{#1}}
\newcommand{\textexponent}[2]{$#1^{\mathrm{#2}}$}

\newcommand{\treeprojection}[1]{p[#1]}

\newcommand{\Tzero}{T_0}
\newcommandx{\union}[2][1 = undefined, 2 = undefined]{
  \ifthenelse{\equal{#1}{undefined}}{\cup}{
    \ifthenelse{\equal{#2}{undefined}}{\bigcup #1}{\bigcup_{#1} #2}
  }
}
\newcommandx{\verticalconcatenation}[2][1 =, 2 =]{
  \ifthenelse{\equal{#1}{}}{\oplus}{\bigoplus_{#1} #2}
}
\newcommand{\verticalsection}[2]{#1_{#2}}

\newcommand{\X}[1][]{\ifthenelse{\equal{#1}{}}{\mathbb{X}}{\mathbb{X}^{#1}}}

\newcommand{\Baire}{Baire\xspace}
\newcommand{\Banach}{Ban\-ach\xspace}

\newcommand{\Borel}{Bor\-el\xspace}
\newcommand{\Conley}{Con\-ley\xspace}
\newcommand{\Dougherty}{Dough\-er\-ty\xspace}

\newcommand{\Fubini}{Fu\-bi\-ni\xspace}

\newcommand{\Harrington}{Har\-ring\-ton\xspace}

\newcommand{\Hjorth}{Hjorth\xspace}
\newcommand{\Jackson}{Jack\-son\xspace}
\newcommand{\Jankov}{Jan\-kov\xspace}
\newcommand{\Kakutani}{Kak\-u\-tan\-i\xspace}
\newcommand{\Kechris}{Kech\-ris\xspace}
\newcommand{\Kuratowski}{Kur\-at\-ow\-ski\xspace}
\newcommand{\Lebesgue}{Leb\-es\-gue\xspace}
\newcommand{\Lecomte}{Le\-comte\xspace}
\newcommand{\Louveau}{Louv\-eau\xspace}
\newcommand{\Lusin}{Lu\-sin\xspace}
\newcommand{\Mazurkiewicz}{Maz\-ur\-kie\-wicz\xspace}
\newcommand{\Miller}{Mil\-ler\xspace}
\newcommand{\Moschovakis}{Mos\-cho\-vak\-is\xspace}
\newcommand{\Mycielski}{My\-ciel\-ski\xspace}
\newcommand{\Novikov}{No\-vik\-ov\xspace}
\newcommand{\Polish}{Po\-lish\xspace}
\newcommand{\Sierpinski}{Sier\-pi\'{n}\-ski\xspace}

\newcommand{\Solecki}{So\-leck\-i\xspace}
\newcommand{\Souslin}{Sous\-lin\xspace}
\newcommand{\Todorcevic}{To\-dor\-cev\-ic\xspace}
\newcommand{\Ulam}{Ul\-am\xspace}
\newcommand{\Vitali}{Vit\-al\i\xspace}
\newcommand{\vonNeumann}{von Neu\-mann\xspace}

\newenvironment{lemmaproof}{
   
  \begin{proof}
}{\end{proof}}

\newenvironment{propositionproof}{
   
  \begin{proof}
}{\end{proof}}

\newenvironment{sublemmaproof}{
   
  \begin{proof}
}{\end{proof}}

\newenvironment{theoremproof}{
   
  \begin{proof}
}{\end{proof}}

\newtheorem{introtheorem}{Theorem}

\newtheorem{conjecture}{Conjecture}[section]
\newtheorem{lemma}[conjecture]{Lemma}
\newtheorem{proposition}[conjecture]{Proposition}
\newtheorem{sublemma}[conjecture]{Sublemma}
\newtheorem{theorem}[conjecture]{Theorem}

\theoremstyle{definition}
\newtheorem*{acknowledgements}{Acknowledgments}

\newtheorem{remark}[conjecture]{Remark}

\begin{document}

\begin{abstract}
  We establish a dichotomy theorem characterizing the circumstances under which a treeable
  \Borel equivalence relation $E$ is essentially countable. Under additional topological 
  assumptions on the treeing, we in fact show that $E$ is essentially countable if and only if there 
  is no continuous embedding of $\Eone$ into $E$. Our techniques also yield the first classical 
  proof of the analogous result for hypersmooth equivalence relations, and allow us to show 
  that up to continuous \Kakutani embeddability, there is a minimum \Borel function which is not 
  essentially countable-to-one.
\end{abstract}

\author[J.D. Clemens]{John D. Clemens}

\address{
  John D. Clemens \\
  Institut f\"{u}r mathematische Logik und Grundlagenforschung \\
  Fachbereich Mathematik und Informatik \\
  Universit\"{a}t M\"{u}nster \\
  Einsteinstra{\ss}e 62 \\
  48149 M\"{u}nster \\
  Germany
 }

\email{jclemens@uni-muenster.de}

\urladdr{
  http://wwwmath.uni-muenster.de/u/john.clemens
}

\author[D. Lecomte]{Dominique Lecomte}

\address{
  Dominique Lecomte \\
  Universit\'{e} Paris 6, Institut de Math\'{e}matiques de Jussieu, Projet Analyse Fonctionnelle \\
  Couloir 16 - 26, 4\`{e}me \'{e}tage, Case 247, 4, place Jussieu, 75 252 Paris Cedex 05, France \\
  and Universit\'{e} de Picardie, I.U.T de l'Oise, site de Creil \\
  13, all\'{e}e de la fa\"iencerie, 60 107 Creil, France
}

\email{dominique.lecomte@upmc.fr}

\urladdr{
  https://www.imj-prg.fr/~dominique.lecomte/
}

\author[B.D. Miller]{Benjamin D. Miller}

\address{
  Benjamin D. Miller \\
  Kurt G\"{o}del Research Center for Mathematical Logic \\
  W\"{a}hringer Stra{\ss}e 25 \\
  1090 Vienna \\
  Austria \\
  and Institut f\"{u}r mathematische Logik und Grundlagenforschung \\
  Fachbereich Mathematik und Informatik \\
  Universit\"{a}t M\"{u}nster \\
  Einsteinstra{\ss}e 62 \\
  48149 M\"{u}nster \\
  Germany
 }

\email{glimmeffros@gmail.com}

\urladdr{
  http://wwwmath.uni-muenster.de/u/ben.miller
}

\thanks{The authors were supported in part by SFB Grant 878.}

  

\title[Essential countability and treeability]{Essential countability of treeable equivalence relations}

\maketitle

\section*{Introduction}

\subsection*{Basic notions}
A \definedterm{\Polish space} is a separable topological space admitting a compatible complete 
metric.  A subset of such a space is $\Ksigma$ if it is a countable union of compact sets, 
\definedsymbol{\Fsigma} if it is a countable union of closed sets, \definedsymbol{\Gdelta} if it is a 
countable intersection of open sets, and \definedterm{\Borel} if it is in the $\sigma$-algebra 
generated by the underlying topology. A function between such spaces is \definedterm
{\Borel} if pre-images of open sets are \Borel. Every subset of a \Polish space inherits the \Borel
structure consisting of its intersection with each \Borel subset of the underlying space.

We endow $\N$ with the discrete topology. A subset of a \Polish space is \definedterm{analytic} if it 
is a continuous image of a closed subset of $\Bairespace$. It is not hard to see that every
non-empty analytic set is a continuous image of $\Bairespace$ itself. A set is \definedterm
{co-analytic} if its complement is analytic. Following standard practice, we use \definedsymbol
{\Sigmaclass[1][1]} and \definedsymbol{\Piclass[1][1]} to denote the classes of analytic and 
co-analytic subsets of \Polish spaces.

Suppose that $X$ and $Y$ are \Polish spaces and $R$ and $S$ are relations on $X$ and $Y$. A 
\definedterm{homomorphism} from $R$ to $S$ is a function $\phi \from X \to Y$ sending 
$R$-related sequences to $S$-related sequences, a \definedterm{cohomomorphism} from $R$ to 
$S$ is a function $\phi \from X \to Y$ sending $R$-unrelated sequences to $S$-unrelated 
sequences, a \definedterm{reduction} of $R$ to $S$ is a function satisying both requirements, and 
an \definedterm{embedding} of $R$ into $S$ is an injective reduction. Given sequences 
$\sequence{R_i}[i \in I]$ and $\sequence{S_i}[i \in I]$ of relations on $X$ and $Y$, we use the 
analogous terminology to describe functions $\phi \from X \to Y$ which have the desired property 
for all $i \in I$.

When $E$ and $F$ are equivalence relations on $X$ and $Y$, the existence of a reduction of
$E$ to $F$ is trivially equivalent to the existence of an injection of $X / E$ into $Y / F$. By requiring
that the former is \Borel, we obtain a definable refinement of cardinality capable of distinguishing
quotients whose classical cardinality is that of $\R$. Over the last few decades, this
notion has been used to great effect in shedding new light on obstacles of definability inherent
in classification problems throughout mathematics, particularly in the theories of countable groups
and fields, probability measure-preserving group actions, separable $C^*$ and von Neumann 
algebras, and separable \Banach spaces. In order to better understand such results, it is essential 
to obtain the best possible understanding of the \Borel reducibility hierarchy. The present paper is
a contribution towards this goal.

\subsection*{An initial segment}
It is easy to see that for each countable cardinal $n$, there is a \Borel reduction of the equality 
relation on the $n$-point discrete space to any \Borel equivalence relation with at least $n$ 
classes. The first non-trivial theorem in the area appears in [Sil80], implying that there is a 
\Borel reduction of the equality relation on $\R$ to any \Borel equivalence relation with 
uncountably many classes. That is, the continuum hypothesis holds in the definable context. 
Building on this, [HKL90, Theorem 1] yields the analog of the continuum 
hypothesis at the next level, namely, that there is a \Borel reduction of the \definedterm{\Vitali 
equivalence relation} on $\R$, i.e., the orbit equivalence relation induced by the action of $\Q$ on 
$\R$ under addition, into any \Borel equivalence relation which is not \Borel reducible to the 
equality relation on $\R$.

Going one step further, [KL97, Theorem 1] implies that under \Borel reducibility,
there is no \Borel equivalence relation lying strictly between the \Vitali equivalence relation and
the orbit equivalence relation induced by the action of $\functions{<\N}{\R}$ on $\functions{\N}
{\R}$ under addition. It is well-known, however, that the full analog of [HKL90, Theorem 1] cannot hold. This can be seen, for example, by noting
that under \Borel reducibility, the latter equivalence relation is incomparable with the orbit
equivalence relation induced by the action of $\functions{\N}{\Q}$ on $\functions{\N}{\R}$.
Nevertheless, in this paper we establish a generalization of [KL97, Theorem 1] of a 
substantially less local nature.

One should note that to facilitate both the proofs of these results as well as topological 
strengthenings in which \Borel reducibility is replaced with continuous embeddability, one
typically focuses on different equivalence relations. In [Sil80], one uses the equality relation
on $\Cantorspace$. In [HKL90], one uses the relation \definedsymbol
{\Ezero} on $\Cantorspace$ given by $x \mathrel{\Ezero} y \iff \exists n \in \N \forall m \ge n
\ x(m) = y(m)$. And in [KL97], one uses the relation \definedsymbol{\Eone} on
$\CantorCantorspace$ given by $x \mathrel{\Eone} y \iff \exists n \in \N \forall m \ge n \ x(m) = 
y(m)$.

\subsection*{Treeable equivalence relations}

We identify graphs with their (ordered) edge sets, so that a \definedterm{graph} on $X$ is an 
irreflexive, symmetric binary relation $G$ on $X$. A \definedterm{cycle} through such a graph is
a sequence $\sequence{x_i}[i \le n]$ such that $n \ge 3$, $\sequence{x_i}[i < n]$ is injective, $x_i 
\mathrel{G} x_{i+1}$ for all $i < n$, and $x_0 = x_n$. We say that $G$ is \definedterm{acyclic} if it 
admits no such cycles. A \definedterm{treeing} of an equivalence relation is an acyclic \Borel 
graph whose connected components coincide with the classes of the relation. A \Borel 
equivalence relation is \definedterm{treeable} if it admits a \Borel treeing. Examples include orbit 
equivalence relations associated with free \Borel actions of countable discrete free groups. Such 
relations play a particularly significant role in the measure-theoretic context, due primarily to their 
susceptability to cocycle reduction techniques.

Beyond such applications, treeable equivalence relations play another important role as a 
proving ground for natural conjectures, where simpler arguments can often be used to obtain
stronger results. One example appears in [Hjo08], in which a strengthening of [HKL90, Theorem 1] is established for treeable \Borel equivalence relations. Although the proof given there takes [HKL90, Theorem 1] for granted, more direct arguments have since appeared (see, for example, [Mil12, Theorem 22]). Moreover, in the presence of strong determinacy assumptions, the ideas behind this argument can be used to establish the natural generalizations of both [Sil80] and [HKL90, Theorem 1] to treeable equivalence relations of higher complexity.

Following the standard abuse of language, we say that an equivalence relation is \definedterm
{finite} if all of its classes are finite, and \definedterm{countable} if all of its classes are countable.
A \Borel equivalence relation is \definedterm{essentially $\calE$} if it is \Borel reducible to a
\Borel equivalence relation in $\calE$. In addition to the results just mentioned, [Hjo08] 
concludes with a question at the heart of our concerns here: is $\Eone$ the minimum treeable 
\Borel equivalence relation which is not essentially countable?

\subsection*{Essential countability}

In order to present our characterization of essential countability, we must first introduce some
terminology. A \definedterm{path} through a binary relation $R$ on $X$ is a sequence of the form 
$\sequence{x_i}[i \le n]$, where $n \in \N$ and $x_i \mathrel{R} x_{i+1}$ for all $i < n$. The 
\definedterm{\textexponent{n}{th} iterate} of $R$ is the binary relation \definedsymbol
{\relationpower{R}{n}} consisting of all pairs $\pair{y}{z}$ for which there is such a path with 
$x_0 = y$ and $x_n = z$. We use \definedsymbol{\relationpower{R}{\le n}} to denote $\union
[m \le n][\relationpower{R}{m}]$.

For all $n \in \N$, let $\F[n]$ denote the equivalence relation on $\CantorCantorspace$ given
by $x \mathrel{\F[n]} y \iff \forall m \ge n \ x(m) = y(m)$.

\begin{introtheorem} \label{introtheorem:homomorphismtreeabledichotomy}
  Suppose that $X$ is a \Polish space, $E$ is a treeable \Borel equivalence relation on $X$, and 
  $G$ is a \Borel treeing of $E$. Then exactly one of the following holds:
  \begin{enumerate}
    \item The equivalence relation $E$ is essentially countable.
    \item There exists a function $f \from \N \to \N$ for which there is a continuous homomorphism
      $\phi \from \CantorCantorspace \to X$ from $\sequence{\F[n+1] \setminus \F[n]}[n \in \N]$ to 
      $\sequence{\relationpower{G}{\le f(n+1)} \setminus \relationpower{G}{\le f(n)}}[n \in \N]$.
  \end{enumerate}
\end{introtheorem}

Although this stops somewhat short of yielding an answer to [Hjo08, Question 13], it does
imply one of the main corollaries of a positive answer: among essentially treeable \Borel
equivalence relations, essential countability is \definedterm{robust}, in the sense that it is 
equivalent to the existence of a universally measurable (or $\aleph_0$-universally \Baire 
measurable) reduction of $E$ to a countable equivalence relation.

Moreover, under appropriate topological assumptions, we do obtain a positive answer to the
original question. We say that a \Borel equivalence relation is \definedterm
{subtreeable-with-$\Fsigma$-iterates} if it has a \Borel treeing contained in an acyclic graph with 
$\Fsigma$ iterates.

\begin{introtheorem} \label{introtheorem:treeabledichotomy}
  Suppose that $X$ is a \Polish space and $E$ is a \Borel equivalence relation on $X$ which is
  essentially subtreeable-with-$\Fsigma$-iterates. Then exactly one of the following holds:
  \begin{enumerate}
    \item The equivalence relation $E$ is essentially countable.
    \item There is a continuous embedding $\pi \from \CantorCantorspace \to X$ of $\Eone$ into 
      $E$.
  \end{enumerate}
\end{introtheorem}

Although the restriction that $E$ is subtreeable-with-$\Fsigma$-iterates might appear rather 
Machiavellian, it turns out that the family of such relations has unbounded potential complexity.
While this fact is beyond the scope of the present paper, one should note that it has surprising 
consequences for the global structure of the \Borel reducibility hierarchy. We say that a class 
$\calE$ of \Borel equivalence relations is \definedterm{dichotomous} if there is a \Borel 
equivalence relation $E_\calE$ such that for every \Borel equivalence relation $E$, either $E \in 
\calE$ or there is a \Borel reduction of $E_\calE$ to $E$. Using the unbounded potential 
complexity of the family of \Borel equivalence relations which are subtreeable-with-$\Fsigma$ 
iterates, one can show that if a \Borel equivalence relation is not \Borel reducible to $\Ezero$, then 
it is incomparable with \Borel equivalence relations of unbounded potential complexity, 
strengthening [KL97, Theorem 2]. It follows that if $\calE$ is a dichotomous class of 
equivalence relations of bounded potential complexity, then $\calE$ consists solely of smooth 
equivalence relations. Consequently, the only non-trivial such families are those associated with 
the main results of [Sil80] and [HKL90], the classes of potentially 
open and potentially closed equivalence relations. These developments will be explored in a 
future paper.

\subsection*{Essentially countable-to-one functions}
In the case of treeings induced by \Borel functions, we obtain even stronger results. To describe
these, we must again introduce some terminology. A \definedterm{\Kakutani embedding} of a 
function $T \from X \to X$ into a function $U \from Y \to Y$ is a \Borel injection $\pi \from X \to Y$ 
with the property that $(\pi \composition T)(x) = (U^n \composition \pi)(x)$, where $n > 0$ is least
with $(U^n \composition \pi)(x) \in \image{\pi}{X}$, for all $x \in X$.

We say that a set $Y \subseteq X$ is \definedterm{$T$-complete} if $X = \union[n \in \N][\image
{T^{-n}}{Y}]$, we say that a set $Y \subseteq X$ is \definedterm{$T$-stable} if $\image{T}{Y} 
\subseteq Y$, and we say that a \Borel function $T \from X \to X$ is \definedterm{essentially 
countable-to-one} if there is a $T$-complete, $T$-stable \Borel set $B \subseteq X$ on which $T$ 
is countable-to-one.

The \definedterm{product} of functions $f \from X \to X$ and $g \from Y \to Y$ is the function $f 
\times g \from X \times Y \to X \times Y$ given by $(f \times g)(x, y) = \pair{f(x)}{g(y)}$. The 
\definedterm{successor function} on $\N$ is given by $S(n) = n + 1$. The \definedterm
{unilateral shift} on $\CantorCantorspace$ is given by $s(\sequence{x_n}[n \in \N]) = \sequence
{x_{n+1}}[n \in \N]$.

\begin{introtheorem} \label{introtheorem:Kakutanidichotomy}
  Suppose that $X$ is a \Polish space and $T \from X \to X$ is a \Borel function. Then exactly one
  of the following holds:
  \begin{enumerate}
    \item The function $T$ is essentially countable-to-one.
    \item There is a continuous \Kakutani embedding $\pi \from \N \times \CantorCantorspace \to
      X$ of $S \times s$ into $T$.
  \end{enumerate}
\end{introtheorem}

\subsection*{Organization}
In \S\ref{section:preliminaries}, we review the basic descriptive set theory utilized throughout. In
\S\ref{section:Baire}, we establish several \Baire category results. In \S\ref{section:Silver}, we 
prove a parametrized version of an unpublished generalization of the main theorem of [Sil80], due originally to \Conley-\Lecomte-\Miller. In \S\ref{section:maindichotomy}, we establish 
our main technical theorems, from which the results stated thus far are relatively straightforward 
corollaries. These technical results are essentially generalizations of Theorems \ref
{introtheorem:homomorphismtreeabledichotomy} and \ref{introtheorem:treeabledichotomy} to 
\Borel equivalence relations equipped with suitably definable assignments of quasi-metrics to 
their classes, although we state them in a somewhat different form so as to facilitate the exposition. 
In \S\ref{section:hypersmooth}, we give the promised classical proof of [KL97, Theorem 1]. In \S\ref{section:treeable}, we establish Theorems \ref
{introtheorem:homomorphismtreeabledichotomy} and \ref{introtheorem:treeabledichotomy}. In
\S\ref{section:Kakutani}, we establish Theorem \ref{introtheorem:Kakutanidichotomy}.

\section{Preliminaries} \label{section:preliminaries}

In this section, we review the basic descriptive set theory utilized throughout the paper. 

Suppose that $X$ and $Y$ are topological spaces. The \definedterm{compact-open topology} on 
the set of all continuous functions $f \from X \to Y$ is that generated by the sets $\set{f \from X \to Y}
[\image{f}{K} \subseteq U]$, where $K \subseteq X$ is compact and $U \subseteq Y$ is open. We 
use \definedsymbol{\continuousfunctions{X}{Y}} to denote the corresponding topological space. 
The following observations will aid complexity calculations involving this space.

\begin{proposition} \label{proposition:evaluation}
  Suppose that $X$ is a compact \Polish space and $Y$ is a \Polish space. Then the function $e 
  \from \continuousfunctions{X}{Y} \times X \to Y$ given by $e(f, x) = f(x)$ is continuous.
\end{proposition}

\begin{propositionproof}
  See, for example, [Kur68, {\S}IV.44.II].
\end{propositionproof}

\begin{proposition} \label{proposition:compactopen}
  Suppose that $X$ is a locally compact \Polish space and $Y$ is a \Polish space. Then
  $\continuousfunctions{X}{Y}$ is a \Polish space.
\end{proposition}

\begin{propositionproof}
  See, for example, [Kur68, {\S}IV.44.VII].
\end{propositionproof}

Although \Borel functions constitute a much broader class than continuous ones, the following 
fact often allows one to treat \Borel functions as if they are continuous.

\begin{proposition} \label{proposition:changeoftopology}
  Suppose that $X$ and $Y$ are \Polish spaces and $\calF$ is a countable family of \Borel
  functions $T \from X \to Y$. Then there are finer \Polish topologies on $X$ and $Y$, 
  whose \Borel sets are the same as those of the original topologies, with respect to which every $T 
  \in \calF$ is continuous. Moreover, if $X = Y$ then the topologies on $X$ and $Y$ can
  be taken to be the same.
\end{proposition}

\begin{propositionproof}
  See, for example, [Kec95, \S13].
\end{propositionproof}

When proving facts about \Polish spaces, it is often notationally convenient (and perhaps 
conceptually clearer) to first focus on the special case of $\Bairespace$. The desired general
result is then typically obtained from a representation theorem such as the following.

\begin{proposition} \label{proposition:surjection}
  Every \Polish space is analytic.
\end{proposition}

\begin{propositionproof}
  See, for example, [Kec95, Theorem 7.9].
\end{propositionproof}

A \definedterm{\Baire space} is a topological space in which every countable intersection of 
dense open sets is dense. The following fact ensures that \Baire category techniques are 
applicable in arbitrary \Polish spaces.

\begin{theorem}[\Baire] \label{theorem:Bairecategorytheorem}
  Every complete metric space is a \Baire space.
\end{theorem}

\begin{theoremproof}
  See, for example, [Kec95, Theorem 8.4].
\end{theoremproof}

A set is \definedterm{nowhere dense} if its closure does not contain a non-empty open set, a set
is \definedterm{meager} if it is contained in a countable union of nowhere dense sets, a set is 
\definedterm{comeager} if its complement is meager, and a set has the \definedterm{\Baire 
property} if its symmetric difference with some open set is meager. One can view the latter three 
notions as topological analogs of $\mu$-null sets, $\mu$-conull sets, and $\mu$-measurable sets, 
although the topological and measure-theoretic notions behave quite differently.

The following fact, known in some circles as \definedterm{localization}, can be viewed as the 
\Baire category analog of the \Lebesgue density theorem.

\begin{proposition} \label{proposition:localization}
  Suppose that $X$ is a \Polish space and $B \subseteq X$ has the \Baire property. Then $B$ is
  non-meager if and only if there is a non-empty open set $U \subseteq X$ such that $B$ is
  comeager in $U$.
\end{proposition}

\begin{propositionproof}
  This easily follows from the definitions of a \Baire space and the \Baire property (see, for 
  example, [Kec95, Proposition 8.26]).
\end{propositionproof}

A function $\phi \from X \to Y$ is \definedterm{\Baire measurable} if for all open sets $V \subseteq 
Y$, the set $\preimage{\phi}{V}$ has the \Baire property. Again, this can be viewed as a
topological analog of $\mu$-measurability. The following observation is a very strong analog 
of the measure-theoretic fact that $\mu$-measurable functions can be approximated by 
continuous ones on sets of arbitrarily large $\mu$-measure.

\begin{proposition} \label{proposition:continuous}
  Suppose that $X$ and $Y$ are \Polish spaces and $\phi \from X \to Y$ is \Baire measurable.
  Then there is a dense $\Gdelta$ set $C \subseteq X$ such that $\restriction{\phi}{C}$ is
  continuous.
\end{proposition}

\begin{propositionproof}
  See, for example, [Kec95, Proposition 8.38].
\end{propositionproof}

Although our primary focus is on \Borel sets, we will often consider analytic sets, in which case the 
following fact ensures that \Baire category arguments remain applicable.

\begin{proposition}[\Lusin-\Sierpinski] \label{proposition:analyticBaire}
  Suppose that $X$ is a \Polish space and $A \subseteq X$ is analytic. Then $A$ has the \Baire
  property.
\end{proposition}

\begin{propositionproof}
  See, for example, [Kec95, Theorem 21.6].
\end{propositionproof}

A topological space $X$ is \definedsymbol{\Tzero} if for all distinct $x, y \in X$, there is an open 
set $U \subseteq X$ containing exactly one of $x$ and $y$. A set $Y \subseteq X$ is \definedterm
{invariant} with respect to an equivalence relation $E$ on $X$ if it is a union of $E$-classes. An 
equivalence relation $E$ on $X$ is \definedterm{generically ergodic} if every invariant set $B
\subseteq X$ with the \Baire property is meager or comeager. The following consequence of 
generic ergodicity is often useful when dealing with param\-etrized dichotomy theorems.

\begin{proposition} \label{proposition:meagertoonehomomorphism}
  Suppose that $X$ is a \Baire space, $Y$ is a second countable $\Tzero$ space, $E$ is a 
  generically ergodic equivalence relation on $X$, and $\phi \from X \to Y$ is a \Baire measurable 
  homomorphism from $E$ to the equality relation on $Y$. Then there exists $y \in Y$ for 
  which $\preimage{\phi}{y}$ is comeager.
\end{proposition}

\begin{propositionproof}
  Fix a basis $\set{V_n}[n \in \N]$ for the topology of $Y$, let $N$ denote the set of $n \in \N$ for
  which $\preimage{\phi}{V_n}$ is comeager, and let $y$ be the unique element of $\intersection
  [n \in N][V_n] \setminus \union[n \in \setcomplement{N}][V_n]$.
\end{propositionproof}

The \definedterm{\textexponent{x}{th} vertical section} and \definedterm{\textexponent{y}{th}
horizontal section} of a set $R \subseteq X \times Y$ are the sets \definedsymbol{\verticalsection
{R}{x}} and \definedsymbol{\horizontalsection{R}{y}} given by $\verticalsection{R}{x} = \set{y \in Y}
[x \mathrel{R} y]$ and $\horizontalsection{R}{y} = \set{x \in X}[x \mathrel{R} y]$. Given a property 
$P$, we write \definedsymbol{\forcomeagerlymany x \ P(x)} to indicate that the set $\set{x \in X}
[P(x)]$ is comeager. The following fact can be viewed as the \Baire category analog of \Fubini's 
theorem.

\begin{theorem}[\Kuratowski-\Ulam] \label{theorem:KuratowskiUlam}
  Suppose that $X$ and $Y$ are \Baire spaces, $Y$ is second countable, and $R \subseteq X 
  \times Y$ has the \Baire property.
  \begin{enumerate}
    \item $\forcomeagerlymany x \in X \ \verticalsection{R}{x}$ has the \Baire property.
    \item $R$ is comeager $\iff \forcomeagerlymany x \in X \ \verticalsection{R}{x}$ is comeager.
  \end{enumerate} 
\end{theorem}

\begin{theoremproof}
  See, for example, [Kec95, Theorem 8.41].
\end{theoremproof}

The following fact can often be used to reduce problems of finding perfect sets with desirable
properties to questions of \Baire category.

\begin{theorem}[\Mycielski] \label{theorem:Mycielski}
  Suppose that $X$ is a non-empty \Polish space and $R \subseteq X \times X$ is meager. Then 
  there is a continuous cohomomorphism $\phi \from \Cantorspace \to X$ from the equality relation
  on $\Cantorspace$ to $R$.
\end{theorem}

\begin{theoremproof}
  See, for example, [Kec95, Theorem 19.1].
\end{theoremproof}

We use \definedsymbol{\sigmaclass{\Sigmaclass[1][1]}} to denote the class of subsets of 
\Polish spaces which lie in the smallest $\sigma$-algebra containing the analytic sets, and we say 
that a function $f \from X \to Y$ is \definedterm{$\sigmaclass{\Sigmaclass[1][1]}$-measurable} if for 
all open sets $U \subseteq Y$, the set $\preimage{f}{U}$ is in $\sigmaclass{\Sigmaclass[1][1]}$. 
We use \definedsymbol{\projection{X}} to denote the \definedterm{projection function} given
by $\projection{X}(x, y) = x$. A \definedterm{uniformization} of a set $R \subseteq X \times Y$ is a 
function $f \from \image{\projection{X}}{R} \to Y$ whose graph is contained in $R$.

\begin{theorem}[\Jankov-\vonNeumann] \label{theorem:JankovvonNeumann}
  Suppose that $X$ and $Y$ are \Polish spaces and $R \subseteq X \times Y$ is an
  analytic set. Then there is a $\sigmaclass{\Sigmaclass[1][1]}$-measurable uniformization of $R$.
\end{theorem}

\begin{theoremproof}
  See, for example, [Kec95, Theorem 18.1].
\end{theoremproof}

\begin{theorem}[\Lusin-\Novikov] \label{theorem:LusinNovikov}
  Suppose that $X$ and $Y$ are \Polish spaces and $R \subseteq X \times Y$ is a \Borel set all of 
  whose vertical sections are countable. Then there are countably many \Borel uniformizations of
  $R$ whose graphs cover $R$.
\end{theorem}

\begin{propositionproof}
  See, for example, [Kec95, Theorem 18.10].
\end{propositionproof}

The following facts will be useful in ensuring that various constructions yield \Borel sets.

\begin{theorem}[\Lusin] \label{theorem:setofunicity}
  Suppose that $X$ and $Y$ are \Polish spaces and $R \subseteq X \times Y$ is \Borel. Then 
  $\set{x \in X}[\exists! y \in Y \ x \mathrel{R} y]$ is co-analytic.
\end{theorem}

\begin{theoremproof}
  See, for example, [Kec95, Theorem 18.11].
\end{theoremproof}

\begin{theorem}[\Lusin] \label{theorem:countabletoone}
  Suppose that $X$ and $Y$ are \Polish spaces and $f \from X \to Y$ is a countable-to-one
  \Borel function. Then $\image{f}{X}$ is \Borel.
\end{theorem}

\begin{theoremproof}
  See, for example, [Kec95, Lemma 18.12].
\end{theoremproof}

Although the class of analytic sets is clearly closed under projections, one must often consider 
analogs of projections in which the non-emptiness of the sections is replaced with stronger 
conditions. The following two facts ensure that the class of analytic sets is also closed under 
certain generalized projections of this form.

\begin{theorem}[\Mazurkiewicz-\Sierpinski] \label{theorem:uncountablequantifier}
  Suppose that $X$ and $Y$ are \Polish spaces and $R \subseteq X \times Y$ is analytic. 
  Then so too is $\set{x \in X}[\verticalsection{R}{x} \text{ is uncountable}]$.
\end{theorem}

\begin{propositionproof}
  See, for example, [Kec95, Theorem 29.20].
\end{propositionproof}

\begin{theorem}[\Novikov] \label{theorem:comeagerquantifier}
  Suppose that $X$ and $Y$ are \Polish spaces and $R \subseteq X \times Y$ is analytic. Then so 
  too is $\set{x \in X}[\verticalsection{R}{x} \text{ is comeager}]$.
\end{theorem}

\begin{theoremproof}
  See, for example, [Kec95, Theorem 29.22].
\end{theoremproof}

Suppose that $\Gamma$ and $\Gamma'$ are classes of subsets of \Polish spaces. A property 
$P$ is \definedterm{$\Gamma$-on-$\Gamma'$} if $\set{x \in X}[P(\verticalsection{R}{x})] \in 
\Gamma$ whenever $X$ and $Y$ are \Polish spaces and $R \subseteq X \times Y$ in $\Gamma'$. 
The following \definedterm{reflection theorem} will help us to ensure that our constructions yield 
\Borel sets.

\begin{theorem}[\Harrington-\Kechris-\Moschovakis] \label{theorem:reflection}
  Suppose that $P$ is a $\Piclass[1][1]$-on-$\Sigmaclass[1][1]$ property. Then every analytic 
  subset of a \Polish space satisfying $P$ is contained in a \Borel set satisfying $P$.
\end{theorem}

\begin{theoremproof}
  See, for example, [Kec95, Theorem 35.10].
\end{theoremproof}

Rather than apply reflection directly, we will often use the following \definedterm{separation 
theorem}.

\begin{theorem}[\Lusin] \label{theorem:separation}
  Suppose that $X$ is a \Polish space and $A, A' \subseteq X$ are disjoint analytic sets. Then there
  is a \Borel set $B \subseteq X$ such that $A \subseteq B$ and $A' \intersection B = \emptyset$.
\end{theorem}

\begin{theoremproof}
  This easily follows from Theorem \ref{theorem:reflection} (a direct proof can be found, for 
  example, in [Kec95, Theorem 14.7]).
\end{theoremproof}

This yields the following connection between analytic and \Borel sets.

\begin{theorem}[\Souslin] \label{theorem:bianalytic}
  A subset of a \Polish space is \Borel if and only if it is both analytic and co-analytic.
\end{theorem}

\begin{theoremproof}
  The fact that sets which are both analytic and co-analytic are \Borel is a direct consequence of 
  Theorem \ref{theorem:separation}, and the converse follows from Propositon \ref
  {proposition:surjection} and a straightforward induction (see, for example, [Kec95, Theorem 14.11], although the latter part is proven there in a somewhat different 
  fashion).
\end{theoremproof}

We will also use the following \definedterm{generalized separation theorem}.

\begin{theorem}[\Novikov] \label{theorem:sequenceseparation}
  Suppose that $X$ is a \Polish space and $A_n \subseteq X$ are analytic sets for which 
  $\intersection[n \in \N][A_n] = \emptyset$. Then there are \Borel sets $B_n \subseteq X$ 
  containing $A_n$ for which $\intersection[n \in \N][B_n] = \emptyset$.
\end{theorem}

\begin{theoremproof}
  This also follows easily from Theorem \ref{theorem:reflection} (a direct proof can be found, for 
  example, in [Kec95, Theorem 28.5]).
\end{theoremproof}

Given $m, n \in \N \union \singleton{\N}$, we say that a sequence $s \in \Cantorspace[m]$ is
\definedterm{extended by} a sequence $t \in \Cantorspace[n]$, or \definedsymbol{s \extendedby
t}, if $s(i) = t(i)$ for all $i < m$. We use \definedsymbol{s \concatenation t} to denote the 
\definedterm{concatenation} of $s$ and $t$.

Fix sequences $s_n \in \Cantorspace[n]$ for which the set $\set{s_n}[n \in \N]$ is \definedterm
{dense}, in the sense that $\forall s \in \Cantortree \exists n \in \N \ s \extendedby s_n$. Let 
\definedsymbol{\Gzero} denote the graph on $\Cantorspace$ consisting of all pairs of the form 
$\pair{s_n \concatenation \singletonsequence{i} \concatenation x}{s_n \concatenation 
\singletonsequence{1-i} \concatenation x}$, where $i < 2$, $n \in \N$, and $x \in \Cantorspace$.

The \definedterm{restriction} of a graph $G$ on $X$ to a set $Y \subseteq X$ is the graph
\definedsymbol{\restriction{G}{Y}} on $Y$ given by $\restriction{G}{Y} = G \intersection (Y \times 
Y)$. Given a graph $G$ on $X$, we say that a set $Y \subseteq X$ is \definedterm
{$G$-independent} if $\restriction{G}{Y} = \emptyset$.

\begin{proposition}[\Kechris-\Solecki-\Todorcevic] \label{proposition:G0independent}
  Suppose that $B \subseteq \Cantorspace$ is a $\Gzero$-independent set with the \Baire 
  property. Then $B$ is meager.
\end{proposition}

\begin{propositionproof}
  This is a direct consequence of the definition of $\Gzero$ and Proposition \ref
  {proposition:localization} (see, for example, [KST99, Proposition 6.2]).
\end{propositionproof}

An \definedterm{$I$-coloring} of $G$ is a function $c \from X \to I$ such that $\preimage{c}
{\singleton{i}}$ is $G$-independent for all $i \in I$. We say that $G$ has \definedterm{countable 
\Borel chromatic number} if there is a \Borel $\N$-coloring of $G$.

\begin{theorem}[\Kechris-\Solecki-\Todorcevic] \label{theorem:G0dichotomy}
  Suppose that $X$ is a \Polish space and $G$ is an analytic graph on $X$. Then exactly one of
  the following holds:
  \begin{enumerate}
    \item The graph $G$ has countable \Borel chromatic number.
    \item There is a continuous homomorphism from $\Gzero$ to $G$.
  \end{enumerate}
\end{theorem}

\begin{theoremproof}
  See, for example, [KST99, Theorem 6.4].
\end{theoremproof}

We say that a \Borel equivalence relation is \definedterm{smooth} if it is \Borel reducible to the 
equality relation on a \Polish space. 

\begin{theorem}[\Harrington-\Kechris-\Louveau] \label{theorem:E0dichotomy}
  Suppose that $X$ is a \Polish space and $E$ is a \Borel equivalence relation on $X$. Then 
  exactly one of the following holds:
  \begin{enumerate}
    \item The equivalence relation $E$ is smooth.
    \item There is a continuous embedding $\pi \from \Cantorspace \to X$ of $\Ezero$ into $E$.
  \end{enumerate}
\end{theorem}

\begin{theoremproof}
  See, for example, [HKL90, Theorem 1.1].
\end{theoremproof}

We say that an equivalence relation is \definedterm{hyper $\calE$} if it is the union of an 
increasing sequence $\sequence{E_n}[n \in \N]$ of relations in $\calE$.

\begin{theorem}[\Dougherty-\Jackson-\Kechris] \label{theorem:countablehypersmooth}
  Suppose that $X$ is a \Polish space and $E$ is a countable \Borel equivalence relation 
  on $X$. If $E$ is hypersmooth, then $E$ is hyperfinite.
\end{theorem}

\begin{theoremproof}
  See, for example, [DJK94, Theorem 5.1].
\end{theoremproof}

We say that a set $B \subseteq X$ is \definedterm{$E$-complete} if it intersects every $E$-class. 
While not strictly necessary for our purposes here, the following fact is also useful in establishing 
closure properties of essential countability.

\begin{theorem}[\Hjorth] \label{theorem:treeableessentiallycountable}
  Suppose that $X$ is a \Polish space and $E$ is a treeable \Borel equivalence relation on $X$.
  Then the following are equivalent:
  \begin{enumerate}
    \item There is an $E$-complete \Borel set on which $E$ is countable.
    \item The equivalence relation $E$ is essentially countable.
  \end{enumerate}
\end{theorem}

\begin{theoremproof}
  See, for example, [Hjo08, Theorem 6].
\end{theoremproof}

Finally, we note that while the original proofs of Theorems \ref{theorem:G0dichotomy}, \ref
{theorem:E0dichotomy}, and \ref{theorem:treeableessentiallycountable} utilized the effective 
theory, classical proofs have since appeared (see [Mil12]). In particular, our reliance on these 
results does not prevent our arguments from being classical in nature.

\section{Baire category results} \label{section:Baire}

In this section, we establish several \Baire category results which will be useful throughout the 
paper.

A function is \definedterm{meager-to-one} if pre-images of singletons are meager.

\begin{proposition} \label{proposition:injection}
  Suppose that $X$ is a \Polish space, $A \subseteq X$, $G$ is a graph on $X$, and there is a
  meager-to-one \Baire measurable function $\phi \from \Cantorspace \to X$ for which the set $A' = 
  \preimage{\phi}{A}$ is comeager and the set $G' = \preimage{(\phi \times \phi)}{G}$ is 
  meager. Then there is a continuous injection $\pi \from \Cantorspace \to A$ of $\Cantorspace$ 
  into a $G$-independent set.
\end{proposition}
  
\begin{propositionproof}
  By Proposition \ref{proposition:continuous}, there is a dense $\Gdelta$ set $B' \subseteq A'$ on 
  which $\phi$ is continuous. Let $E'$ denote the pullback of the equality relation on $X$ through 
  $\phi$. The fact that $\phi$ is \Baire measurable ensures that $E'$ has the \Baire property, and 
  the fact that $\phi$ is meager-to-one implies that every vertical section of $E'$ is meager, so $E'$ 
  is meager by Theorem \ref{theorem:KuratowskiUlam}. In particular, it follows that $(B' \times B') 
  \setminus (E' \union G')$ is a comeager subset of $\Cantorspace \times \Cantorspace$, so 
  Theorem \ref{theorem:Mycielski} yields a continuous injection $\psi \from \Cantorspace \to B'$ of 
  $\Cantorspace$ into a $G'$-independent set which is also a \definedterm{partial transversal} of 
  $E'$, meaning that it intersects every equivalence class in at most one point. It follows that the 
  function $\pi = \phi \composition \psi$ is as desired.
\end{propositionproof}

Throughout the paper, we will work with spaces of the form $\CantorCantorspace[m][n]$, where
$m, n \in \N \union \singleton{\N}$. We use \definedsymbol{\horizontalconcatenation} to denote 
horizontal concatenation, and \definedsymbol{\verticalconcatenation} to denote vertical 
concatenation. We will abuse language by saying that a sequence $s \in \CantorCantorspace[m]
[n]$ is \definedterm{extended} by a sequence $s' \in \CantorCantorspace[m'][n']$, or $s 
\extendedby s'$, if $\forall i < m \forall j < n \ s(i)(j) = s'(i)(j)$.
  
\begin{proposition} \label{proposition:meagerk}
  Suppose that $k \in \N$ and $B \subseteq \CantorCantorspace$ is a set with the \Baire property
  on which $\F[k+1]$ has countable index over $\F[k]$. Then $B$ is meager.
\end{proposition}

\begin{propositionproof}
  Suppose, towards a contradiction, that $B$ is non-meager. Then Theorem \ref
  {theorem:KuratowskiUlam} yields a non-meager set of $\pair{x}{z} \in 
  \CantorCantorspace[k] \times \CantorCantorspace$ such that $\set{y \in \Cantorspace}
  [x \concatenation \singletonsequence{y} \concatenation z \in B]$ is non-meager, and
  therefore uncountable. As $\pair{x \concatenation \singletonsequence{y} \concatenation z}
  {x \concatenation \singletonsequence{y'} \concatenation z} \in \F[k+1] \setminus \F[k]$ for
  distinct $y, y' \in \Cantorspace$, this contradicts the fact that $\F[k+1]$ has countable index over 
  $\F[k]$ on $B$.
\end{propositionproof}

\begin{remark}
  Suppose that $\mu$ is a \Borel probability measure on $\CantorCantorspace$ for which 
  $\mu$-almost every measure in the disintegration of $\mu$ with respect to the function 
  deleting the \textexponent{k}{th} column is continuous (this holds, for example, if $\mu(U) = 1 / 
  2^n$ for every basic open set $U \subseteq \CantorCantorspace$ specifying values on $n$ 
  coordinates). Then an essentially identical argument (using this assumption in place of Theorem
  \ref{theorem:KuratowskiUlam}) yields the analogous result in which $B$ is $\mu$-measurable 
  instead of \Baire measurable.
\end{remark}

\begin{proposition} \label{proposition:meager}
  Suppose that $A \subseteq \CantorCantorspace$ is an analytic set on which $\Eone$ has 
  countable index over $\F[k]$, for some $k \in \N$. Then $\saturation{A}{\Eone}$ is meager.
\end{proposition}

\begin{propositionproof}
  Note that $\F[\ell+1]$ has countable index over $\F[\ell]$ on $\saturation{A}{\F[\ell]}$, for all $\ell 
  \ge k$. As each $\saturation{A}{\F[\ell]}$ is analytic, Proposition \ref{proposition:analyticBaire} 
  ensures that it has the \Baire property, so Proposition \ref{proposition:meagerk} implies that it is 
  meager, thus so too is the set $\saturation{A}{\Eone} = \union[\ell \ge k][\saturation{A}{\F[\ell]}]$.
\end{propositionproof}

Suppose that $X$ is a \Polish space and $E$ is a \Borel equivalence relation on $X$. Theorem
\ref{theorem:LusinNovikov} immediately implies that if $B \subseteq X$ is an $E$-complete \Borel
set on which $E$ is countable, then there is a \Borel reduction of $E$ to $\restriction{E}{B}$, thus
$E$ is essentially countable. Together with Proposition \ref{proposition:meager}, the following 
weak converse yields a simple proof of [KL97, Proposition 1.4], ruling out the 
existence of a \Baire measurable reduction of $\Eone$ to a countable equivalence relation on a 
\Polish space.

\begin{proposition} \label{proposition:essentiallycountable}
  Suppose that $X$ and $Y$ are \Polish spaces, $E$ is an analytic equivalence relation on $X$, 
  $F$ is a countable equivalence relation on $Y$, and there is a \Baire measurable reduction 
  $\phi \from X \to Y$ of $E$ to $F$. Then there is a \Borel set $B \subseteq X$ such that 
  $\restriction{E}{B}$ is countable and $\saturation{B}{E}$ is comeager.
\end{proposition}

\begin{propositionproof}
  By Proposition \ref{proposition:continuous}, there is a dense $\Gdelta$ set $C \subseteq X$ on 
  which $\phi$ is continuous. By Theorem \ref{theorem:JankovvonNeumann}, there is a 
  $\sigmaclass{\Sigmaclass[1][1]}$-measurable function $\phi' \from \image{\phi}{C} \to C$ such 
  that $\phi \composition \phi'$ is the identity function. As pre-images of analytic sets under
  continuous functions are analytic, it follows that $\phi' \composition \phi$ is also
  $\sigmaclass{\Sigmaclass[1][1]}$-measurable, so one more application of Proposition \ref
  {proposition:continuous} yields a dense $\Gdelta$ set $D \subseteq C$ on which it is continuous. 
  Then the set $A = \image{(\phi' \composition \phi)}{D}$ is analytic. As $E$ is countable on $A$ 
  and Theorem \ref{theorem:uncountablequantifier} ensures that this property is $\Piclass[1]
  [1]$-on-$\Sigmaclass[1][1]$, Theorem \ref{theorem:reflection} yields a \Borel set $B \supseteq 
  A$ on which $E$ is countable. As $D \subseteq \saturation{B}{E}$, it follows that the latter set is 
  comeager.
\end{propositionproof}

For each $k \in \N$, let \definedsymbol{\F[k][m][n]} denote the equivalence relation on
$\CantorCantorspace[m][n]$ given by $x \mathrel{\F[k][m][n]} y \iff \forall i \ge k \ x(i) = y(i)$. We 
say that $\phi \from \CantorCantorspace[m][n] \to \CantorCantorspace[m'][n']$ is \definedterm
{extended by} $\psi \from \CantorCantorspace[m''][n''] \to \CantorCantorspace[m'''][n''']$, or $\phi 
\extendedby \psi$, if $s \extendedby t \implies \phi(s) \extendedby \psi(t)$ for all $s \in 
\CantorCantorspace[m][n]$ and $t \in \CantorCantorspace[m''][n'']$.

\begin{proposition} \label{proposition:onedimensionalMycielskigeneralization}
  Suppose that $m, m', n \in \N$, $\phi \from \CantorCantorspace[m][n] \to \CantorCantorspace[m']
  [n]$ is an embedding of $\sequence{\F[k][m][n]}[k < n]$ into $\sequence{\F[k][m'][n]}[k < n]$, and 
  $\calU$ is a family of open subsets of $\CantorCantorspace[n]$ whose union is dense. Then 
  there exists $m'' \in \N$ for which there is an embedding $\psi \from \CantorCantorspace[m]
  [n] \to \CantorCantorspace[m''][n]$ of $\sequence{\F[k][m][n]}[k < n]$ into $\sequence{\F[k][m''][n]}
  [k < n]$ extending $\phi$ with $\forall s \in \CantorCantorspace[m][n] \exists U \in \calU
  \ \extensions{\psi(s)} \subseteq U$.
\end{proposition}

\begin{propositionproof}
  Fix an injective enumeration $\sequence{s_i}[i < I]$ of $\CantorCantorspace[m][n]$. Set $m_0 = 
  m'$ and $\phi_0 = \phi$, and recursively find $m_{i+1} \in \N$ and $\phi_{i+1} \from 
  \CantorCantorspace[m][n] \to \CantorCantorspace[m_{i+1}][n]$ of the form $\phi_{i+1}(s) = \phi_i
  (s) \verticalconcatenation t$, where $t \in \CantorCantorspace[m_{i+1} - m_i][n]$ has the property 
  that $\extensions{\phi_i(s_i) \verticalconcatenation t}$ is a subset of some $U \in \calU$. Set $m'' 
  = m_I$ and $\psi = \phi_I$.
\end{propositionproof}

\begin{proposition} \label{proposition:onedimensionalconfigurationstabilization}
  Suppose that $m, n \in \N$ and $\pi \from \CantorCantorspace[n] \to \N$ is \Baire measurable. 
  Then there exist $i \from \CantorCantorspace[m][n] \to \N$, $m' \in \N$, and an embedding $\phi 
  \from \CantorCantorspace[m][n] \to \CantorCantorspace[m'][n]$ of $\sequence{\F[k][m][n]}[k < n]$ 
  into $\sequence{\F[k][m'][n]}[k < n]$, extending the identity function on $\CantorCantorspace[m]
  [n]$, with the property that $\forall s \in \CantorCantorspace[m][n] \forcomeagerlymany x \in 
  \CantorCantorspace[n] \ i(s) = \pi(\phi(s) \verticalconcatenation x)$.
\end{proposition}

\begin{propositionproof}
  Proposition \ref{proposition:localization} ensures that the family $\calU$ of open sets $U
  \subseteq \CantorCantorspace[n]$ with the property that $\exists i \in \N \forcomeagerlymany x \in 
  U \ i = \pi(x)$ has dense union. The desired result therefore follows from an application of 
  Proposition \ref{proposition:onedimensionalMycielskigeneralization} to the identity function on 
  $\CantorCantorspace[m][n]$.
\end{propositionproof}

Note that $\F[k+1][n] \setminus \F[k][n]$ is homeomorphic to the product of $(\CantorCantorspace[k]
[n] \times \CantorCantorspace[k][n]) \times \CantorCantorspace[n-(k+1)]$ with the complement of
the equality relation on $\Cantorspace$. In particular, it is a locally compact \Polish space, so
Theorem \ref{theorem:Bairecategorytheorem} ensures that it is a \Baire space.

\begin{proposition} \label{proposition:twodimensionalMycielskigeneralization}
  Suppose that $\ell, m, m', n \in \N$, $\phi \from \CantorCantorspace[m][n] \to \CantorCantorspace
  [m'][n]$ is an embedding of $\sequence{\F[k][m][n]}[k < n]$ into $\sequence{\F[k][m'][n]}[k < n]$, 
  and $\calU$ is a family of open subsets of $\F[\ell+1][n]$ whose union is dense. Then there exists 
  $m'' \in \N$ for which there is an embedding $\psi \from \CantorCantorspace[m][n] \to 
  \CantorCantorspace[m''][n]$ of $\sequence{\F[k][m][n]}[k < n]$ into $\sequence{\F[k][m''][n]}[k < 
  n]$ extending $\phi$ with the property that $\forall \pair{s}{t} \in \F[\ell+1][m][n] \setminus \F[\ell][m]
  [n] \exists U \in \calU \ \extensions{\psi(s)} \times \extensions{\psi(t)} \subseteq U$.
\end{proposition}

\begin{propositionproof}
  Fix an injective enumeration $\sequence{s_i, t_i}[i < I]$ of $\F[\ell+1][m][n] \setminus \F[\ell][m][n]$. 
  Define $m_0 = m'$ and $\phi_0 = \phi$, and recursively find $m_{i+1} \in \N$ and $\phi_{i
  +1} \from \CantorCantorspace[m][n] \to \CantorCantorspace[m_{i+1}][n]$ of the form $\phi_{i+1}
  (s) = \phi_i(s) \verticalconcatenation \sigma(s)$, where $\sigma \from \CantorCantorspace[m][n] 
  \to \CantorCantorspace[m_{i+1} - m_i][n]$ is itself of the form
  \begin{equation*}
    \sigma(s) =
      \begin{cases}
        t & \text{if $s \mathrel{\F[\ell][m][n]} s_i$ and } \\
        u & \text{otherwise,}
      \end{cases}
  \end{equation*}
  and $\F[\ell+1][n] \intersection (\extensions{\phi_i(s_i) \verticalconcatenation t} 
  \times \extensions{\phi_i(t_i) \verticalconcatenation u})$ is a non-empty subset of some $U \in 
  \calU$. Set $m'' = m_I$ and $\psi = \phi_I$.
\end{propositionproof}

\begin{proposition} \label{proposition:twodimensionalconfigurationstabilization}
  Suppose that $\ell, m, n \in \N$ and 
  $$\pi \from \F[\ell+1][n] \setminus \F[\ell][n] \to \N$$ is \Baire 
  measurable. Then there exist $i \from \F[\ell+1][m][n] \setminus \F[\ell][m][n] \to \N$, $m' \in \N$, 
  and an embedding $\phi \from \CantorCantorspace[m][n] \to \CantorCantorspace[m'][n]$ of 
  $\sequence{\F[k][m][n]}[k < n]$ into $\sequence{\F[k][m'][n]}[k < n]$, extending the identity function
  on $\CantorCantorspace[m][n]$, with the property that\smallskip 
  
  \leftline{$\forall \pair{s}{t} \in \F[\ell+1][m][n] 
  \setminus \F[\ell][m][n] \forcomeagerlymany \pair{x}{y} \in \F[\ell+1][n]$}\smallskip
  
  \rightline{$i(s, t) = \pi(\phi(s) 
  \verticalconcatenation x, \phi(t) \verticalconcatenation y)$.}
\end{proposition}

\begin{propositionproof}
  By Proposition \ref{proposition:localization}, the family $\calU$ of open sets $U 
  \subseteq \F[\ell+1][n]$ with the property that $\exists i \in \N 
  \forcomeagerlymany \pair{x}{y} \in U \ i = \pi(x, y)$ has dense union. The desired result therefore 
  follows from an application of Proposition \ref
  {proposition:twodimensionalMycielskigeneralization} to the identity function on 
  $\CantorCantorspace[m][n]$.
\end{propositionproof}

We next establish an analog of Theorem \ref{theorem:Mycielski} for $\CantorCantorspace[n]$.

\begin{proposition} \label{proposition:longgeneralizedMycielski}
  Suppose that $m, m', n \in \N$, $\phi \from \CantorCantorspace[m][n] \to \CantorCantorspace[m']
  [n]$, $C \subseteq \CantorCantorspace[n]$ is comeager, and $\sequence{R_k}[k < n]$ is a 
  sequence of subsets of $\CantorCantorspace[n] \times \CantorCantorspace[n]$ with the property 
  that $R_k$ is comeager in $\F[k+1][n]$, for all $k < n$. Then $\phi$ extends to a continuous 
  homomorphism $\psi \from \CantorCantorspace[n] \to C$ from $\sequence{\F[k][n], \F[k+1][n] 
  \setminus \F[k][n]}[k < n]$ to $\sequence{\F[k][n], R_k}[k < n]$.
\end{proposition}

\begin{propositionproof}
  Fix a sequence $\sequence{U_i}[i \in \N]$ of dense open subsets of $\CantorCantorspace[n]$ 
  whose intersection is contained in $C$. For all $k < n$, fix a decreasing sequence $\sequence
  {U_{i,k}}[i \in \N]$ of dense open subsets of $\F[k+1][n]$ whose intersection is contained in $R_k$. 
  We will recursively construct a strictly increasing sequence of natural numbers $\ell_i \in \N$ and 
  embeddings $\phi_i \from \CantorCantorspace[i][n] \to \CantorCantorspace[\ell_i][n]$ of 
  $\sequence{\F[k][i][n]}[k < n]$ into $\sequence{\F[k][\ell_i][n]}[k < n]$ with the following properties:
  \begin{enumerate}
    \item $\forall s \in \CantorCantorspace[i][n] \forall t \in \CantorCantorspace[i+1][n] 
      \ (s \extendedby t \implies \phi_i(s) \extendedby \phi_{i+1}(t))$.
    \item $\forall s \in \CantorCantorspace[i+1][n] \ \extensions{\phi_{i+1}(s)} \subseteq U_i$.
    \item $\forall k < n \forall \pair{s}{t} \in \F[k+1][i+1][n] \setminus \F[k][i+1][n] \\
      \hspace*{7pt} \F[k+1][n] \intersection (\extensions{\phi_{i+1}(s)} \times \extensions{\phi_{i+1}(t)}) 
        \subseteq U_{i,k}$.
  \end{enumerate}
  We begin by setting $\ell_m = m'$ and $\phi_m = \phi$. Given $\phi_i \from \CantorCantorspace[i]
  [n] \to \CantorCantorspace[m_i][n]$, define $\phi_i' \from \CantorCantorspace[i+1][n] \to 
  \CantorCantorspace[m_i + 1][n]$ by $\phi_i'(s \verticalconcatenation t) = \phi_i(s) 
  \verticalconcatenation t$. We then obtain $m_{i+1} \in \N$ and $\phi_{i+1} \from 
  \CantorCantorspace[i+1][n] \to \CantorCantorspace[m_{i+1}][n]$ by one application of 
  Proposition \ref{proposition:onedimensionalMycielskigeneralization} and $n$ applications of 
  Proposition \ref{proposition:twodimensionalMycielskigeneralization}. This completes the 
  recursive construction, and the corresponding function $\psi \from \CantorCantorspace[n] \to 
  \CantorCantorspace[n]$, given by $\psi(x) = \union[i \ge m][\phi_i \composition \projection
  {\CantorCantorspace[i][n]}(x)]$, is as desired.
\end{propositionproof}

We next consider analogous results with $\CantorCantorspace$ in place of $\CantorCantorspace
[n]$.

\begin{proposition} \label{proposition:longonedimensionalMycielskigeneralization}
  Suppose that $m, m', n, n' \in \N$, $\phi \from \CantorCantorspace[m][n] \to \CantorCantorspace
  [m'][n']$ is an embedding of $\sequence{\F[k][m][n]}[k \le n]$ into $\sequence{\F[k][m'][n']}[k \le n]$, 
  and $\calU$ is a family of open subsets of $\CantorCantorspace$ whose union is dense. Then 
  there exist $m'', n'' \in \N$ for which there is an embedding $\psi \from \CantorCantorspace[m]
  [n] \to \CantorCantorspace[m''][n'']$ of $\sequence{\F[k][m][n]}[k \le n]$ into $\sequence{\F[k][m'']
  [n'']}[k \le n]$ extending $\phi$ with the property that $\forall s \in \CantorCantorspace[m][n] \exists 
  U \in \calU \ \extensions{\psi(s)} \subseteq U$.
\end{proposition}

\begin{propositionproof}
  Fix an injective enumeration $\sequence{s_i}[i < I]$ of $\CantorCantorspace[m][n]$. Set $m_0 = 
  m'$, $n_0 = n'$, and $\phi_0 = \phi$, and recursively find $m_{i+1}, n_{i+1} \in \N$ and 
  $\phi_{i+1} \from \CantorCantorspace[m][n] \to \CantorCantorspace[m_{i+1}][n_{i+1}]$ of the form 
  $\phi_{i+1}(s) = (\phi_i(s) \horizontalconcatenation u) \verticalconcatenation v$, where $u \in 
  \CantorCantorspace[m_i][n_{i+1} - n_i]$ and $v \in \CantorCantorspace[m_{i+1} - m_i][n_{i+1}]$
  have the property that $\extensions{(\phi_i(s_i) \horizontalconcatenation u) \verticalconcatenation 
  v}$ is a subset of some $U \in \calU$. Set $m'' = m_I$, $n'' = n_I$, and $\psi = \phi_I$.
\end{propositionproof}

We say that an open set $U \subseteq \CantorCantorspace \times \CantorCantorspace$ is 
\definedterm{$k$-dense} if for all $m, n \in \N$ and $\pair{s}{t} \in \setcomplement{\F[k][m][n]}$, 
there exist $m', n' \in \N$ and extensions $s', t' \in \CantorCantorspace[m'][n']$ of $s, t$ such that
$\extensions{s'} \times \extensions{t'} \subseteq U$ and
\begin{equation*}
  \forall i < m' \forall k < j < n' \ (s'(j)(i) \neq t'(j)(i) \implies (i < m \mathand j < n)).
\end{equation*}

\begin{proposition} \label{proposition:longtwodimensionalMycielskigeneralization}
  Suppose that $\ell, m, m', n, n' \in \N$, $$\phi \from \CantorCantorspace[m][n] \to 
  \CantorCantorspace[m'][n']$$ is an embedding of $\sequence{\F[k][m][n]}[k \le n]$ into $\sequence
  {\F[k][m'][n']}[k \le n]$, and $\calU$ is a family of open subsets of $\CantorCantorspace \times 
  \CantorCantorspace$ whose union is $\ell$-dense. Then there exist $m'', n'' \in \N$ for which 
  there is an embedding $\psi \from \CantorCantorspace[m][n] \to \CantorCantorspace[m''][n'']$ of 
  $\sequence{\F[k][m][n]}[k \le n]$ into $\sequence{\F[k][m''][n'']}[k \le n]$ extending $\phi$ with the 
  property that $\forall \pair{s}{t} \in \setcomplement{\F[\ell][m][n]} \exists U \in \calU \ \extensions{\psi
  (s)} \times \extensions{\psi(t)} \subseteq U$.
\end{proposition}

\begin{propositionproof}
  Fix an injective enumeration $\sequence{s_i, t_i}[i < I]$ of $\setcomplement{\F[\ell][m][n]}$. Define 
  $m_0 = m'$, $n_0 = n'$, and $\phi_0 = \phi$, and recursively find $m_{i+1}, n_{i+1}
  \in \N$ and $\phi_{i+1} \from \CantorCantorspace[m][n] \to \CantorCantorspace[m_{i+1}]
  [n_{i+1}]$ of the form $\phi_{i+1}(s) = (\phi_i(s) \horizontalconcatenation \sigma(s))
  \verticalconcatenation \tau(s)$, where $\sigma \from \CantorCantorspace[m][n] \to
  \CantorCantorspace[m_i][n_{i+1} - n_i]$ is of the form
  \begin{equation*}
    \sigma(s) =
      \begin{cases}
        t & \text{if $s \mathrel{\F[\ell][m][n]} s_i$ and } \\
        u & \text{otherwise,}
      \end{cases}
  \end{equation*}
  $\tau \from \CantorCantorspace[m][n] \to \CantorCantorspace[m_{i+1} - m_i][n_{i+1}]$ is of
  the form  
  \begin{equation*}
    \tau(s) =
      \begin{cases}
        v & \text{if $s \mathrel{\F[\ell][m][n]} s_i$ and } \\
        w & \text{otherwise,}
      \end{cases}
  \end{equation*}  
  and $\pair{t}{u} \in \F[\max(0, \ell+1-n_i)][m_i][n_{i+1} - n_i]$ and 
  $$\pair{v}{w} \in \F[\ell+1][m_{i+1} - 
  m_i][n_{i+1}]$$ are such that $\extensions{(\phi_i(s_i) \horizontalconcatenation t) 
  \verticalconcatenation v} \times \extensions{(\phi_i(t_i) \horizontalconcatenation u) 
  \verticalconcatenation w}$ is contained in some $U \in \calU$. Set $m'' = m_I$ and $\psi = \phi_I$.
\end{propositionproof}

We say that a set $M \subseteq \CantorCantorspace$ is \definedterm{$k$-meager} if it is
disjoint from the intersection of a countable family of $k$-dense open sets.

\begin{proposition} \label{proposition:longlonggeneralizedMycielski}
  Suppose that $C \subseteq \CantorCantorspace$ is comeager and $\sequence{R_k}[k \in \N]$ is 
  a sequence of subsets of $\CantorCantorspace \times \CantorCantorspace$ with the property 
  that $R_k$ is $k$-meager, for all $k \in \N$. Then there is a continuous homomorphism $\phi \from 
  \CantorCantorspace \to C$ from $\sequence{\F[k], \setcomplement{\F[k]}}[k \in \N]$ to $\sequence
  {\F[k], \setcomplement{R_k}}[k \in \N]$.
\end{proposition}

\begin{propositionproof}
  Fix a sequence $\sequence{U_i}[i \in \N]$ of dense open subsets of $\CantorCantorspace$ 
  whose intersection is contained in $C$. For all $k \in \N$, fix a decreasing sequence $\sequence
  {U_{i,k}}[i \in \N]$ of $k$-dense open subsets of $\CantorCantorspace \times 
  \CantorCantorspace$ whose intersection is disjoint from $R_k$. We will recursively construct
  strictly increasing sequences of natural numbers $m_i, n_i \in \N$ and embeddings $\phi_i \from 
  \CantorCantorspace[i][i] \to \CantorCantorspace[m_i][n_i]$ of $\sequence{\F[k][i][i]}[k \le i]$ into 
  $\sequence{\F[k][m_i][n_i]}[k \le i]$ such that:
  \begin{enumerate}
    \item $\forall s \in \CantorCantorspace[i][i] \forall t \in \CantorCantorspace[i+1][i+1] 
      \ (s \extendedby t \implies \phi_i(s) \extendedby \phi_{i+1}(t))$.
    \item $\forall s \in \CantorCantorspace[i+1][i+1] \ \extensions{\phi_{i+1}(s)} \subseteq U_i$.
    \item $\forall k \le i \forall \pair{s}{t} \in \setcomplement{\F[k][i+1][i+1]} \ \extensions{\phi_{i+1}(s)} 
      \times \extensions{\phi_{i+1}(t)} \subseteq U_{i,k}$.
  \end{enumerate}
  We begin by setting $m_0 = n_0 = 0$ and fixing $\phi_0 \from \CantorCantorspace[0][0] \to 
  \CantorCantorspace[0][0]$. Given $\phi_i \from \CantorCantorspace[i][i] \to \CantorCantorspace
  [m_i][n_i]$, define $\psi_i \from \CantorCantorspace[i+1][i+1] \to \CantorCantorspace[m_i + 1]
  [n_i + 1]$ by $\psi_i((s \horizontalconcatenation t) \verticalconcatenation u) = (\phi_i(s)
  \horizontalconcatenation t) \verticalconcatenation u$. We then obtain $m_{i+1}, n_{i+1} \in \N$ 
  and $\phi_{i+1} \from \CantorCantorspace[i+1][i+1] \to \CantorCantorspace[m_{i+1}][n_{i+1}]$ by 
  one application of Proposition \ref{proposition:longonedimensionalMycielskigeneralization} and 
  $i+1$ applications of Proposition \ref{proposition:longtwodimensionalMycielskigeneralization}. 
  This completes the recursive construction. Define $\phi \from \CantorCantorspace \to 
  \CantorCantorspace$ by $\phi(x) = \union[i \in \N][\phi_i \composition \projection
  {\CantorCantorspace[i][i]}(x)]$.
\end{propositionproof}

We next give a condition sufficient for ensuring $k$-meagerness.

\begin{proposition} \label{proposition:containmentmeager}
  Suppose that $k \in \N$ and $R \subseteq \CantorCantorspace \times \CantorCantorspace$ is an 
  $\Fsigma$ set disjoint from $\Eone \setminus \F[k]$. Then $R$ is $k$-meager.
\end{proposition}

\begin{propositionproof}
  It is sufficient to show that every open set $U \subseteq \Cantorspace \times \Cantorspace$
  containing $\Eone \setminus \F[k]$ is $k$-dense. Towards this end, suppose that $m, n \in 
  \N$ and $\pair{s}{t} \in \setcomplement{\F[k][m][n]}$. Let $x, y \in \CantorCantorspace$ denote
  the extensions of $s, t$ with constant value $0$ off of the domains of $s, t$. Then $\pair{x}{y} \in
  \Eone \setminus \F[k]$, so $\pair{x}{y} \in U$, thus there exist $m', n' \in \N$ and $s', t' \in 
  \CantorCantorspace[m'][n']$ such that $s \extendedby s' \extendedby x$, $t \extendedby t' 
  \extendedby y$, and $\extensions{s'} \times \extensions{t'} \subseteq U$.
\end{propositionproof}

We close this section with a closure property of the family of equivalence relations into which
$\Eone$ is reducible.

\begin{proposition} \label{proposition:Eoneclosure}
  Suppose that $X$ and $Y$ are \Polish spaces, $E$ and $F$ are \Borel equivalence relations on
  $X$ and $Y$, $A \subseteq X$ is analytic, and $\phi \from A \to Y$ is a \Borel reduction of $E$ to 
  $F$ for which there is a \Baire measurable reduction $\psi \from \CantorCantorspace \to \image
  {\phi}{A}$ of $\Eone$ to $\restriction{F}{\image{\phi}{A}}$. Then there is a continuous embedding 
  of $\Eone$ into $\restriction{E}{A}$.
\end{proposition}

\begin{propositionproof}
  By Proposition \ref{proposition:continuous}, there is a dense $\Gdelta$ set $C \subseteq 
  \CantorCantorspace$ on which $\psi$ is continuous. By Theorem \ref
  {theorem:JankovvonNeumann}, there is a $\sigmaclass{\Sigmaclass[1][1]}$-measurable function 
  $\phi' \from \image{\phi}{A} \to X$ for which $\phi \composition \phi'$ is the identity function. Then 
  $\phi' \composition (\restriction{\psi}{C})$ is a $\sigmaclass{\Sigmaclass[1][1]}$-measurable
  reduction of $\restriction{\Eone}{C}$ to $E$. One more appeal to Proposition \ref
  {proposition:continuous} therefore yields a dense $\Gdelta$ set $D \subseteq C$ for which
  it is a continuous reduction of $\restriction{\Eone}{D}$ to $E$. As Propositions \ref
  {proposition:longlonggeneralizedMycielski} and \ref{proposition:containmentmeager} ensure 
  that there is a continuous embedding of $\Eone$ into $\restriction{\Eone}{D}$, the proposition 
  follows.
\end{propositionproof}

\section{Independent perfect sets} \label{section:Silver}

We say that a set $B \subseteq Y$ is \definedterm{$\aleph_0$-universally \Baire} if $\preimage{f}
{B}$ has the \Baire property whenever $X$ is a \Polish space and $f \from X \to Y$ is continuous. 
In this section, we establish a local version of the following generalization of the perfect set 
theorem for co-analytic equivalence relations (see [Sil80]).

\begin{proposition}[\Conley-\Lecomte-\Miller] \label{proposition:perfectpair}
  Suppose that $X$ is a \Polish space, $A \subseteq X$ is analytic, $G$ is an 
  $\aleph_0$-universally \Baire graph on $X$, $R$ is a reflexive symmetric co-analytic binary 
  relation on $X$, and $\relationpower{G}{2} \subseteq R$. Then at least one of the following 
  holds:
  \begin{enumerate}
    \item There is a \Borel set $B \supseteq A$ on which $\setcomplement{R}$ has countable 
      \Borel chromatic number.
    \item There is a continuous injection $\pi \from \Cantorspace \to A$ of $\Cantorspace$ into a 
      $G$-independent set.
  \end{enumerate}
\end{proposition}

\begin{propositionproof}
  As the property of being independent with respect to an analytic graph is $\Piclass[1]
  [1]$-on-$\Sigmaclass[1][1]$, Theorem \ref{theorem:reflection} ensures that every 
  $(\setcomplement{R})$-independent analytic set is contained in a $(\setcomplement
  {R})$-independent \Borel set. It follows that if $\chromaticnumber[Borel]{\restriction
  {\setcomplement{R}}{A}} \le \aleph_0$, then there is a \Borel set $B \supseteq A$ for which 
  $\chromaticnumber[Borel]{\restriction{\setcomplement{R}}{B}} \le \aleph_0$. 
  Otherwise, Theorem \ref{theorem:G0dichotomy} yields a continuous homomorphism 
  $\phi \from \Cantorspace \to X$ from $\Gzero$ to $(A \times A) \intersection \setcomplement{R}$. 
  As $\Gzero$ has full projection, it follows that $\image{\phi}{\Cantorspace} \subseteq A$, and
  Proposition \ref{proposition:G0independent} ensures that $\phi$ is meager-to-one. So by
  Proposition \ref{proposition:injection}, it only remains to verify that the graph $G' = \preimage{(\phi 
  \times \phi)}{G}$ is meager. Suppose, towards a contradiction, that this is not the case. By 
  Theorem \ref{theorem:KuratowskiUlam}, there exists $x \in X$ for which $\verticalsection{G'}{x}$
  is non-meager and has the \Baire property. Proposition \ref{proposition:G0independent} then 
  yields a pair $\pair{y}{z} \in \restriction{\Gzero}{\verticalsection{G'}{x}}$, in which case the fact that 
  $\relationpower{G}{2} \subseteq R$ implies that $\phi(y) \mathrel{R} \phi(z)$, contradicting the 
  fact that $\phi$ is a homomorphism from $\Gzero$ to $\setcomplement{R}$.
\end{propositionproof}

We will need the following complexity calculation.

\begin{proposition} \label{proposition:complexity}
  Suppose that $X$ is a \Polish space and $G$ is a co-analytic graph on $X$. Then the property 
  $P(A)$ that there is no continuous injection of $\Cantorspace$ into a $G$-independent subset of 
  $A$ is $\Piclass[1][1]$-on-$\Sigmaclass[1][1]$.
\end{proposition}

\begin{propositionproof}
  Let $E$ denote the equality relation on $X$, and suppose that $Y$ is a \Polish space and $R 
  \subseteq X \times Y$ is analytic. Then Proposition \ref{proposition:injection} ensures that the
  inexistence of a continuous injection of $\Cantorspace$ into a $G$-independent subset of 
  $\horizontalsection{R}{y}$ is equivalent to the inexistence of a continuous function $\phi \from 
  \Cantorspace \to X$ for which $\preimage{(\phi \times \phi)}{(\horizontalsection{R}
  {y} \times \horizontalsection{R}{y}) \setminus (E \union G)}$ is comeager. Propositions \ref
  {proposition:evaluation} and \ref{proposition:compactopen} along with Theorem \ref
  {theorem:comeagerquantifier} imply that this latter property is $\Piclass[1][1]$-on-$\Sigmaclass[1]
  [1]$.
\end{propositionproof}

Given an equivalence relation $E$ on $X$, we say that a graph $G$ on $X$ has \definedterm
{countable $E$-local \Borel chromatic number} if its restriction to each equivalence class of $E$ 
has countable \Borel chromatic number.

\begin{proposition} \label{proposition:reflection}
  Suppose that $X$ is a \Polish space, $A \subseteq X$ is analytic, $E$ is an analytic equivalence 
  relation on $X$, $G$ is a co-analytic graph on $X$, $R$ is a reflexive symmetric co-analytic
  binary relation on $X$, and $\relationpower{G}{2} \subseteq R$. Then at least one of the 
  following holds:
  \begin{enumerate}
    \item There is a \Borel set $B \supseteq A$ on which $\setcomplement{R}$ has countable 
      $E$-local \Borel chromatic number.
    \item There exists $x \in X$ for which there is a continuous injection $\pi \from \Cantorspace \to 
      A$ of $\Cantorspace$ into a $G$-independent subset of $\equivalenceclass{x}{E}$.
  \end{enumerate}
\end{proposition}

\begin{propositionproof}
  By Proposition \ref{proposition:complexity}, the property $Q(A)$ that there is no $x \in X$ for 
  which there is a continuous injection $\pi \from \Cantorspace \to A$ of $\Cantorspace$ into a 
  $G$-independent subset of $\equivalenceclass{x}{E}$ is $\Piclass[1][1]$-on-$\Sigmaclass[1][1]$. 
  So if condition (2) fails, then Theorem \ref{theorem:reflection} yields a \Borel set $B \subseteq X$ 
  containing $A$ such that there is no $x \in X$ for which there is a continuous injection $\pi \from 
  \Cantorspace \to B$ of $\Cantorspace$ into a $G$-independent subset of $\equivalenceclass{x}
  {E}$. As Proposition \ref{proposition:analyticBaire} ensures that $G$ is $\aleph_0$-universally 
  \Baire, Proposition \ref{proposition:perfectpair} implies that $\setcomplement{R}$ has countable 
  $E$-local \Borel chromatic number on $B$.
\end{propositionproof}

\section{Two dichotomy theorems} \label{section:maindichotomy}

In this section, we establish the main technical results of the paper. We say that a sequence 
$\sequence{G_n}[n \in \N]$ \definedterm{eventually} has a property $P$ if $G_n$ has property 
$P$ for all but finitely many $n \in \N$.

\begin{theorem} \label{theorem:Boreldichotomy}
  Suppose that $X$ is a \Polish space, $E$ is an analytic equivalence relation on $X$, and 
  $\sequence{R_n}[n \in \N]$ is an increasing sequence of reflexive symmetric co-analytic
  binary relations on $X$ such that $E \subseteq \union[n \in \N][R_n]$ and \heightcorrection
  {$\relationpower{R_n}{2} \subseteq R_{n+1}$} for all $n \in \N$. Then exactly one of the following 
  holds:
  \begin{enumerate}
    \item The set $X$ is a countable union of \Borel sets on which $\sequence
      {\setcomplement{R_n}}[n \in \N]$ eventually has countable $E$-local \Borel chromatic number.
    \item There exists $f \from \N \to \N$ for which there is a continuous homomorphism $\phi \from
      \CantorCantorspace \to X$ from $\sequence{\F[n+1] \setminus \F[n]}[n \in \N]$ to $\sequence{E 
      \intersection R_{f(n+1)} \setminus R_{f(n)}}[n \in \N]$.
  \end{enumerate}
\end{theorem}

\begin{theoremproof}
  Observe that if $f \from \N \to \N$ and $\phi \from \CantorCantorspace \to X$ is a homomorphism
  from $\sequence{\F[n+1] \setminus \F[n]}[n \in \N]$ to $\sequence{E \intersection R_{f(n+1)} 
  \setminus R_{f(n)}}[n \in \N]$, then $\phi$ is necessarily a homomorphism from $\Eone$ to $E$. 
  Moreover, as each of the sets $\F[n+1] \setminus \F[n]$ is non-empty and $\sequence{R_n}[n \in 
  \N]$ is increasing, it follows that $f(n+1) > f(n)$ for all $n \in \N$, so $f(n) \ge n$ for all $n \in \N$, 
  thus $\phi$ is in fact a homomorphism from $\sequence{\Eone \setminus \F[n]}[n \in \N]$ to 
  $\sequence{E \setminus R_n}[n \in \N]$.
  
  To see that conditions (1) and (2) are mutually exclusive, observe that if both hold, then there is 
  a non-meager analytic set $A \subseteq \CantorCantorspace$ such that for all $x \in A$, there 
  is an $\aleph_0$-coloring $c$ of $\restriction{\setcomplement{R_n}}{\equivalenceclass{\phi(x)}
  {\restriction{E}{\image{\phi}{A}}}}$. Then for any such $x$ and $c$, the function $c \composition 
  \phi$ is a coloring of $\restriction{\setcomplement{\F[n]}}{\equivalenceclass{x}{\restriction{\Eone}
  {A}}}$, so $\Eone$ has countable index over $\F[n]$ on $A$, contradicting Proposition \ref
  {proposition:meager}.

  In order to show that at least one of conditions (1) and (2) does indeed hold, it will be convenient 
  to assume that $X = \Bairespace$. To see that this special case is sufficient to establish the 
  theorem, note that we can assume $X$ is non-empty, in which case Proposition \ref
  {proposition:surjection} yields a continuous surjection $\pi \from \Bairespace \to X$. Set $E' = 
  \preimage{(\pi \times \pi)}{E}$ and $R_n' = \preimage{(\pi \times \pi)}{R_n}$. If $f \from \N \to \N$ 
  and $\phi' \from \CantorCantorspace \to \Bairespace$ is a continuous homomorphism from 
  $\sequence{\F[n+1] \setminus \F[n]}[n \in \N]$ to $\sequence{E' \intersection R_{f(n+1)}' \setminus 
  R_{f(n)}'}[n \in \N]$, then the map $\phi = \pi \composition \phi'$ is a homomorphism from 
  $\sequence{\F[n+1] \setminus \F[n]}[n \in \N]$ to $\sequence{E \intersection R_{f(n+1)} \setminus 
  R_{f(n)}}[n \in \N]$. On the other hand, suppose there are \Borel sets $B_n' \subseteq 
  \Bairespace$ and natural numbers $k_n \in \N$ such that $\Bairespace = \union[n \in \N][B_n']$
  and $\restriction{\setcomplement{R_{k_n}'}}{B_n'}$ has countable $E'$-local \Borel chromatic 
  number for all $n \in \N$. Then $X$ is the union of the analytic sets $A_n = \image{\pi}{B_n'}$. If 
  $x \in A_n$, then there exists $x' \in B_n'$ such that $\pi(x') = x$, and if $c' \from B_n' \to \N$ is a 
  coloring of $\restriction{\setcomplement{R_{k_n}'}}{\equivalenceclass{x'}{\restriction{E'}{B_n'}}}$, 
  then the function $c(y) = \min \set{c(y')}[y = \pi(y')]$ is a coloring of $\restriction{\setcomplement
  {R_{k_n}}}{\equivalenceclass{x}{\restriction{E}{A_n}}}$, so Proposition \ref
  {proposition:reflection} yields \Borel sets $B_n \supseteq A_n$ such that $\restriction
  {\setcomplement{R_{k_n+1}}}{B_n}$ has countable $E$-local \Borel chromatic number for all $n
  \in \N$.
    
  We now proceed to the main argument. We will recursively define a decreasing sequence 
  $\sequence{X^\alpha}[\alpha < \omega_1]$ of \Borel subsets of $X$, beginning with $X^0 = X$, 
  and taking intersections at limit ordinals. In order to describe the construction of $X^{\alpha+1}$ 
  from $X^\alpha$, we need several preliminaries.
  
  \begin{lemma}
    There is an increasing sequence $\sequence{S_n}[n \in \N]$ of reflexive symmetric analytic 
    binary relations on $X$ such that $\forall n \in \N \ S_n \subseteq R_n$, $E = \union[n \in \N]
    [S_n]$, and \heightcorrection{$\forall n \in \N \ \relationpower{S_n}{2} \subseteq S_{n+1}$}.
  \end{lemma}
  
  \begin{lemmaproof}
    As $\sequence{E \setminus R_n}[n \in \N]$ is a sequence of analytic sets with empty
    intersection, Theorem \ref{theorem:sequenceseparation} yields a sequence $\sequence{R_n'}
    [n \in \N]$ of \Borel sets with empty intersection such that $\forall n \in \N \ E \setminus R_n 
    \subseteq R_n'$. By replacing $R_n'$ with $\intersection[i \le n][R_i']$, we can ensure that 
    $\sequence{R_n'}[n \in \N]$ is decreasing. By replacing $R_n'$ with $\set{\pair{x}{y} \in R_n'}[x 
    \neq y \mathand \pair{y}{x} \in R_n']$, we can assume that each of these sets is irreflexive and 
    symmetric. Set $R_n'' = E \setminus R_n'$. Then $\sequence{R_n''}[n \in \N]$ is an increasing 
    sequence of reflexive symmetric analytic binary relations on $X$ such that $\forall n \in \N \ R_n'' 
    \subseteq R_n$ and $E = \union[n \in \N][R_n'']$. Set $S_0 = R_0''$, and recursively define 
    \heightcorrection{$S_{n+1} = R_{n+1}'' \union \relationpower{S_n}{2}$}. A straightforward 
    induction shows that $\forall n \in \N \ S_n \subseteq R_n \intersection S_{n+1}$, and it is clear 
    that $E = \union[n \in \N][S_n]$ and \heightcorrection{$\forall n \in \N \ \relationpower{S_n}
    {2} \subseteq S_{n+1}$}.
  \end{lemmaproof}
  
  Fix trees $T_{m,n}$ on $(\N \times \N) \times \N$ for which $\treeprojection
  {T_{m,n}} = S_n \setminus R_m$. An \definedterm{approximation} is a quadruple of the form $a = 
  \quadruple{n}{f}{\phi}{\sequence{\psi_k}[k < n]}$, with the property that $n \in \N$, $f \from \set{0, 
  \ldots, n} \to \N$, $\phi \from \CantorCantorspace[n][n] \to \Bairespace[n]$, and $\psi_k \from \F[k
  +1][n][n] \setminus \F[k][n][n] \to \Bairespace[n]$ for all $k < n$. 
  
  We say that $a$ is \definedterm{extended} by another approximation $b$ if $n^a \le n^b$, $f^a 
  \extendedby f^b$, $\phi^a \extendedby \phi^b$, and $\psi^a \extendedby \psi^b$ for all $k <
  n^a$. When $n^b = n^a + 1$, we say that $b$ is a \definedterm{one-step extension} of $a$.
  
  A \definedterm{configuration} is a quadruple of the form $\gamma = \quadruple{n}{f}{\phi}
  {\sequence{\psi_k}[k < n]}$, with the property that $n \in \N$, $f \from \set{0, \ldots, n} \to \N$, $\phi 
  \from \CantorCantorspace[n] \to \Bairespace$, and $\psi_k \from \F[k+1][n] \setminus \F[k][n] \to 
  \Bairespace$ for all $k < n$. 
  
  For reasons of definability, it will be important to focus our attention on configurations which are 
  \definedterm{continuous}, in the sense that the functions $\phi$ and $\psi_k$ are continuous. In 
  the course of the argument, it will also be useful to consider configurations which are merely 
  \definedterm{\Baire measurable}, in the sense that the functions $\phi$ and $\psi_k$ are \Baire 
  measurable.
  
  We say that $\gamma$ is \definedterm{compatible} with a set $Y \subseteq X$ if 
  $\phi^\gamma(x) \in Y$ for all $x \in \domain{\phi^\gamma}$. We say that $\gamma$ is 
  \definedterm{compatible} with the sequence $\sequence{T_{m,n}}[m,n \in \N]$ if $\pair{\pair
  {\phi^\gamma(x)}{\phi^\gamma(y)}}{\psi_k^\gamma(x,y)} \in \branches{T_{f(k), f(k+1)}}$ for all 
  $k < n^\gamma$ and $\pair{x}{y} \in \domain{\psi_k^\gamma}$. We say that $\gamma$ is 
  \definedterm{compatible} with an approximation $a$ if $n^a = n^\gamma$, $f^a = f^\gamma$,
  $\phi^a \extendedby \phi^\gamma$, and $\psi_k^a \extendedby \psi_k^\gamma$ for all $k <
  n^a$. 
  
  Again for reasons of definability, it will be important to focus on the corresponding notions of 
  \definedterm{generic compatibility}, in which one only asks for the desired properties on a
  comeager set. Although it is possible to proceed with only this latter notion, the arguments 
  provide a strong connection between the two, and only a modicum of further effort is 
  required to elucidate the connection between them.
  
  Given an embedding $\pi \from \CantorCantorspace[n^\gamma] \to \CantorCantorspace[n^
  \gamma]$ of $\sequence{\F[k][n^\gamma]}[k < n^\gamma]$ into $\sequence{\F[k][n^\gamma]}[k < 
  n^\gamma]$, let \definedsymbol{\pushforward{\pi}{\gamma}} denote the configuration $\delta$ 
  given by $n^\delta = n^\gamma$, $f^\delta = f^\gamma$, $\phi^\delta = \phi^\gamma 
  \composition \pi$, and $\psi_k^\delta = \psi_k^\gamma \composition (\pi \times \pi)$.
  
  \begin{lemma} \label{lemma:compatible}
    Suppose that $a$ is an approximation, $B \subseteq X$ is a \Borel set, $\gamma$ is a
    \Baire measurable configuration which is generically compatible with $a$, $B$, and $\sequence
    {T_{m,n}}[m,n \in \N]$, $m, m' \in \N$, and $\pi \from \CantorCantorspace[m][n^\gamma] \to
    \CantorCantorspace[m'][n^\gamma]$ is an embedding of $\sequence{\F[k][m]
    [n^\gamma]}[k < n^\gamma]$ into $\sequence{\F[k][m'][n^\gamma]}[k < n^\gamma]$. Then $\pi$
    extends to a continuous embedding $\pi' \from \CantorCantorspace[n^\gamma] \to 
    \CantorCantorspace[n^\gamma]$ of $\sequence{\F[k][n^\gamma]}[k < n^\gamma]$ into
    $\sequence{\F[k][n^\gamma]}[k < n^\gamma]$ for which $\pushforward{\pi'}{\gamma}$ is 
    continuous and compatible with $a$, $B$, and $\sequence{T_{m,n}}[m,n \in \N]$.
  \end{lemma}
  
  \begin{lemmaproof}
    By Proposition \ref{proposition:continuous}, there are comeager sets $C \subseteq \domain
    {\phi^\gamma}$ and $C_k \subseteq \domain{\psi_k^\gamma}$ for which $\restriction{\phi^
    \gamma}{C}$ and $\restriction{\psi_k^\gamma}{C_k}$ are continuous. Then the set $D = 
    \preimage{(\phi^\gamma)}{B} \intersection C$ is comeager, as are the sets $D_k \subseteq 
    \domain{\psi_k^\gamma}$ of $\pair{x}{y} \in C_k$ with $\pair{\pair{\phi^\gamma(x)}{\phi^\gamma
    (y)}}{\psi_k^\gamma(x, y)} \in \branches{T_{f^\gamma(k), f^\gamma(k+1)}}$ and $\pair{s}{t} \in 
    \domain{\psi_k^a} \implies \psi_k^a(s, t) \extendedby \psi_k^\gamma(x, y)$, where $s$ and $t$ 
    are the projections of $x$ and $y$ onto $\domain{\phi^a}$. But Proposition \ref
    {proposition:longgeneralizedMycielski} ensures that the function $\pi$ extends to a continuous 
    homomorphism $\pi' \from \CantorCantorspace[n^\gamma] \to D$ from $\sequence{\F[k]
    [n^\gamma], \F[k+1][n^\gamma] \setminus \F[k][n^\gamma]}[k < n^\gamma]$ to $\sequence{\F[k]
    [n^\gamma], D_k}[k < n^\gamma]$, and any such function is as desired.
  \end{lemmaproof}

  Given a natural number $n \in \N$ and an embedding $\pi \from \CantorCantorspace[n^\gamma]
  [n^\gamma] \to \CantorCantorspace[n][n^\gamma]$ of $\sequence{\F[k][n^\gamma][n^\gamma]}
  [k < n^\gamma]$ into $\sequence{\F[k][n][n^\gamma]}[k < n^\gamma]$, let \definedsymbol
  {\pushforward{\pi}{\gamma}} denote the configuration $\delta$ given by $n^\delta = n^\gamma$, 
  $f^\delta = f^\gamma$, $\phi^\delta(s \verticalconcatenation x) = \phi^\gamma(\pi(s) 
  \verticalconcatenation x)$, and $\psi_k^\delta(s \verticalconcatenation x, t \verticalconcatenation 
  y) = \psi_k^\gamma(\pi(s) \verticalconcatenation x, \pi(t) \verticalconcatenation y)$.
  
  \begin{lemma} \label{lemma:genericallycompatible}
    Suppose that $\gamma$ is a \Baire measurable configuration. Then there exists $n \in \N$
    for which there is an embedding $\pi \from \CantorCantorspace[n^\gamma][n^\gamma] \to 
    \CantorCantorspace[n][n^\gamma]$ of $\sequence{\F[k][n^\gamma][n^\gamma]}[k < n^\gamma]$ 
    into $\sequence{\F[k][n][n^\gamma]}[k < n^\gamma]$ with the property that $\pushforward{\pi}
    {\gamma}$ is generically compatible with an approximation.
  \end{lemma}
  
  \begin{lemmaproof}
    This follows from one application of Proposition \ref
    {proposition:onedimensionalconfigurationstabilization} and $n^\gamma$ applications of Proposition
    \ref{proposition:twodimensionalconfigurationstabilization}.
  \end{lemmaproof}
    
  Let \definedsymbol{\Gamma^\alpha(a)} denote the set of all continuous configurations which are
  generically compatible with $X^\alpha$, $\sequence{T_{m,n}}[m, n \in \N]$, and $a$. Theorem
  \ref{theorem:comeagerquantifier} ensures that $\Gamma^\alpha(a)$ is analytic (and even \Borel).
  
  Associate with each configuration $\gamma$ the set $D^\gamma \subseteq 
  \CantorCantorspace[n^\gamma]$ given by
  \begin{equation*}
    D^\gamma = \set{x \in \CantorCantorspace[n^\gamma]}[{\forcomeagerlymany y \in 
      \CantorCantorspace[n^\gamma] \ \phi^\gamma(x) \mathrel{S_{f^\gamma(n^\gamma)}} \phi^
        \gamma(y)}].
  \end{equation*}
  If $\gamma$ is generically compatible with $\sequence{T_{m,n}}[m, n \in \N]$, then $D^\gamma$ 
  is comeager. As \heightcorrection{$\relationpower{S_{f^\gamma(n^\gamma)}}{2} \subseteq 
  S_{f^\gamma(n^\gamma)+1}$}, it follows that $\image{\phi^\gamma}{D^\gamma}$ is an 
  $S_{f^\gamma(n^\gamma)+1}$-clique.   
  
  We say that $a$ is \definedterm{$\alpha$-terminal} if $\Gamma^\alpha(b) = \emptyset$ for all 
  one-step extensions $b$ of $a$. Define $A^\alpha(a) = \union[\gamma \in \Gamma^\alpha(a)]
  [\image{\phi^\gamma}{D^\gamma}]$. 

  \begin{lemma} \label{lemma:terminal}
    Suppose that $a$ is an approximation for which there is a continuous injection $\pi \from 
    \Cantorspace \to A^\alpha(a)$ into an $(E \setminus R_{f^a(n^a)+2})$-clique. Then $a$ is not 
    $\alpha$-terminal.
  \end{lemma}
  
  \begin{lemmaproof}
    We first note that $E$ can be replaced with an appropriate $S_n$.
    
    \begin{sublemma}
      There exists $n > f^a(n^a) + 2$ for which there is a continuous injection $\pi' \from 
      \Cantorspace \to \Cantorspace$ such that $\image{(\pi \composition \pi')}{\Cantorspace}$ is an 
      $S_n$-clique.
    \end{sublemma}
    
    \begin{sublemmaproof}
      Fix $x \in X$ with $\image{\pi}{\Cantorspace} \subseteq \equivalenceclass{x}{E}$, and set
      $S_m' = \preimage{\pi}{\horizontalsection{S_m}{x}}$ for all $m \in \N$. Then $\Cantorspace 
      = \union[m \in \N][S_m']$, so there exists $m \ge f^a(n^a) + 2$ for which $S_m'$ is non-meager.
      As Proposition \ref{proposition:analyticBaire} ensures that $S_m'$ has the \Baire property, 
      the one-dimensional analog of Theorem \ref{theorem:Mycielski} (whose proof is even simpler
      than in the two-dimensional case) yields a continuous injection $\pi' \from \Cantorspace \to 
      S_m'$. Set $n = m + 1$. As $\image{(\pi \composition \pi')}{\Cantorspace} \subseteq 
      \horizontalsection{S_m}{x}$ and \heightcorrection{$\relationpower{S_m}{2} \subseteq S_n$}, it 
      follows that $\image{(\pi \composition \pi')}{\Cantorspace}$ is an $S_n$-clique.
    \end{sublemmaproof}
    
    Replacing $\pi$ with $\pi \composition \pi'$, we can assume that $\image{\pi}{\Cantorspace}$ is 
    an $S_n$-clique.
    
    \begin{sublemma}
      There is a continuous injection $\pi_\Gamma \from \Cantorspace \to \Gamma^\alpha(a)$ for 
      which there is a continuous injection $\pi' \from \Cantorspace \to \Cantorspace$ with the 
      property that $(\pi \composition \pi')(x) \in \image{\phi^{\pi_\Gamma(x)}}{D^{\pi_\Gamma(x)}}$
      for all $x \in \Cantorspace$. 
    \end{sublemma}
    
    \begin{sublemmaproof}
      Note that the set of pairs $\pair{x}{\gamma} \in \Cantorspace \times \Gamma^\alpha(a)$ with 
      $\pi(x) \in \image{\phi^\gamma}{D^\gamma}$ is analytic (and even \Borel). By Theorem \ref
      {theorem:JankovvonNeumann}, there is a $\sigmaclass{\Sigmaclass[1][1]}$-measurable 
      function $\pi_\Gamma' \from \Cantorspace \to \Gamma^\alpha(a)$ such that $\pi(x) \in \image
      {\phi^{\pi_\Gamma'(x)}}{D^{\pi_\Gamma'(x)}}$ for all $x \in \Cantorspace$. As $\pi$ is injective 
      and no two distinct points of $\image{\pi}{\Cantorspace}$ are $R_{f^a(n^a)+2}$-related, it 
      follows that $\pi_\Gamma'$ is injective. By Proposition \ref{proposition:continuous}, there is a 
      comeager set $C \subseteq \Cantorspace$ on which $\pi_\Gamma'$ is continuous. The 
      one-dimensional analog of Theorem \ref{theorem:Mycielski}
      therefore yields a continuous injection $\pi' \from \Cantorspace \to C$. Set $\pi_\Gamma = 
      \pi_\Gamma' \composition \pi'$.
    \end{sublemmaproof}
    
    Replacing $\pi$ with $\pi \composition \pi'$, we can assume that $\pi(x) \in \image{\phi^
    {\pi_\Gamma(x)}}{D^{\pi_\Gamma(x)}}$ for all $x \in \Cantorspace$. Note that $\phi^{\pi_\Gamma
    (x)}(x') \mathrel{(S_{n+2} \setminus R_{f^a(n^a)})} \phi^{\pi_\Gamma(y)}(y')$ whenever $x, y \in 
    \Cantorspace$ are distinct, $x' \in D^{\pi_\Gamma(x)}$, and $y' \in D^{\pi_\Gamma(y)}$. Observe 
    further that by 
    Proposition \ref{proposition:evaluation}, the set of pairs $\pair{\quadruple{x}{x'}{y}{y'}}{z} \in 
    (\Cantorspace \times \CantorCantorspace[n^a] \times \Cantorspace \times \CantorCantorspace
    [n^a]) \times \Bairespace$ with the property that $\pair{\pair{\phi^{\pi_\Gamma(x)}(x')}{\phi^
    {\pi_\Gamma(y)}(y')}}{z} \in \branches{T_{f^a(n^a), n+2}}$ is closed, so by Theorem \ref
    {theorem:JankovvonNeumann}, there is a $\sigmaclass{\Sigmaclass[1][1]}$-measurable 
    function $\psi \from \F[n^a+1][n^a+1] \setminus \F[n^a][n^a+1] \to \Bairespace$ such that
    \begin{equation*}
      \pair{\pair{\phi^{\pi_\Gamma(x)}(x')}{\phi^{\pi_\Gamma(y)}(y')}}{\psi(x' \horizontalconcatenation 
        \singletonsequence{x}, y' \horizontalconcatenation \singletonsequence{y})} \in \branches
          {T_{f^a(n^a), n+2}}
    \end{equation*}
    for all distinct $x, y \in \Cantorspace$, $x' \in D^{\pi_\Gamma(x)}$, and $y' \in D^{\pi_\Gamma(y)}$.
    
    Let $\gamma$ denote the \Baire measurable configuration given by $n^\gamma = n^a + 1$, 
    $\restriction{f^\gamma}{\set{0, \ldots, n^a}} = f^a$, $f^\gamma(n^\gamma) = n + 2$, 
    \heightcorrection{$\phi^\gamma(x \horizontalconcatenation \singletonsequence{z}) =
    \phi^{\pi_\Gamma(z)}(x)$}, \heightcorrection{$\psi_k^\gamma(x \horizontalconcatenation
    \singletonsequence{z}, y \horizontalconcatenation \singletonsequence{z}) =
    \psi_k^{\pi_\Gamma(z)}(x, y)$} for $k < n^\gamma$, and $\psi_{n^a}^\gamma = \psi$. Lemma 
    \ref{lemma:genericallycompatible} then yields an approximation $b$, a natural 
    number $n'$, and an embedding $\pi' \from \CantorCantorspace[n^\gamma][n^\gamma] \to 
    \CantorCantorspace[n'][n^\gamma]$ of $\sequence{\F[k][n^\gamma][n^\gamma]}[k < 
    n^\gamma]$ into $\sequence{\F[k][n'][n^\gamma]}[k < n^\gamma]$, extending the identity
    function on $\CantorCantorspace[n^\gamma][n^\gamma]$, for which $\pushforward{\pi'}
    {\gamma}$ is generically compatible with $b$, $\sequence{T_{m,n}}[m,n \in \N]$, and
    $X^\alpha$. As $\pushforward{\pi'}{\gamma}$ is made up of perfectly many configurations
    generically compatible with $a$, it follows that $b$ is a one-step extension of $a$. As 
    Lemma \ref{lemma:compatible} yields a continuous embedding $\pi'' \from 
    \CantorCantorspace[n^\gamma] \to \CantorCantorspace[n^\gamma]$ of $\sequence{\F[k][n^
    \gamma]}[k < n^\gamma]$ into $\sequence{\F[k][n^\gamma]}[k < n^\gamma]$, extending $\pi'$,
    with the property that $\pushforward{\pi''}{\gamma}$ is continuous and compatible with $b$,
    $\sequence{T_{m,n}}[m,n \in \N]$, and $X^\alpha$, it follows that $a$ is not $\alpha$-terminal.
  \end{lemmaproof}
  
  \begin{lemma} \label{lemma:Borelterminal}
    Suppose that $a$ is an $\alpha$-terminal approximation. Then there is a \Borel set $B
    \subseteq X$ containing $A^\alpha(a)$ on which $\setcomplement{R_{f^a(n^a)+3}}$ has
    countable $E$-local \Borel chromatic number.
  \end{lemma}
  
  \begin{lemmaproof}
    As $A^\alpha(a)$ is analytic and Lemma \ref{lemma:terminal} ensures that there is no 
    continuous injection $\pi \from \Cantorspace \to A^\alpha(a)$ into an $(E \setminus R_{f^a(n^a)
    +2})$-clique, the desired result follows from Proposition \ref{proposition:reflection}.
  \end{lemmaproof}
  
  Let \definedsymbol{\terminal{\alpha}} denote the set of all $\alpha$-terminal approximations.
  For every such approximation $a$, appeal to Lemma \ref{lemma:Borelterminal} to obtain a
  \Borel set $B^\alpha(a) \subseteq X$ containing $A^\alpha(a)$ on which $\setcomplement
  {R_{f^a(n^a)+3}}$ has countable $E$-local \Borel chromatic number. Define $X^{\alpha + 1} =
  X^\alpha \setminus \union[a \in \terminal{\alpha}][B^\alpha(a)]$. This completes the recursive 
  construction.
    
  \begin{lemma} \label{lemma:extension}
    Suppose that $a$ is an approximation whose one-step extensions are all $\alpha$-terminal.
    Then $a$ is $(\alpha + 1)$-terminal.  
  \end{lemma}
  
  \begin{lemmaproof}
    Suppose that $b$ is a one-step extension of $a$. If $\gamma$ is a continuous configuration
    generically compatible with $b$, then the $\alpha$-terminality of $b$ ensures that $\image
    {\phi^\gamma}{D^\gamma} \intersection X^{\alpha+1} = \emptyset$. It follows that if $\gamma$
    is also generically compatible with $\sequence{T_{m,n}}[m,n \in \N]$, then it is not generically 
    compatible with $X^{\alpha + 1}$, thus $\Gamma^{\alpha+1}(b) = \emptyset$.
  \end{lemmaproof}
  
  Note that the family of $\alpha$-terminal approximations is increasing. As there are only 
  countably many approximations, there exists $\alpha < \omega_1$ such that every $(\alpha + 
  1)$-terminal approximation is $\alpha$-terminal. If the unique approximation $a$ for which $n^a 
  = f^a(0) = 0$ is $\alpha$-terminal, then $X^{\alpha+1} = \emptyset$, and Lemma \ref
  {lemma:Borelterminal} ensures that there are \Borel sets $B_n \subseteq X$ and natural 
  numbers $k_n \in \N$ such that $X = \union[n \in \N][B_n]$ and $\setcomplement{R_{k_n}}$ has 
  countable $E$-local \Borel chromatic number on $B_n$, for all $n \in \N$.

  Otherwise, Lemma \ref{lemma:extension} allows us to recursively construct 
  non-$\alpha$-term\-inal approximations $a_n$ with the property that $a_{n+1}$ is a one-step 
  extension of $a_n$. Define $f \from \N \to \N$ by $f(n) = f^{a_n}(n)$; define $\phi \from
  \CantorCantorspace \to \Bairespace$ by $\restriction{\phi(x)}{n} = \phi^{a_n}(s)$, where $s$
  is the projection of $x$ onto $\CantorCantorspace[n][n]$; and define $\psi_k \from \F[k+1] 
  \setminus \F[k] \to \Bairespace$ by $\restriction{\psi_k(x, y)}{n} = \psi_k^{a_n}(s, t)$, where $k < n$
  and $s$ and $t$ are the projections of $x$ and $y$ onto $\CantorCantorspace[n][n]$.
  
  It remains to show that for all $k \in \N$, the function $\phi$ is a homomorphism from $\F[k+1] 
  \setminus \F[k]$ to $S_{f(k+1)} \setminus R_{f(k)}$. Towards this end, suppose that $x
  \mathrel{\F[k+1] \setminus \F[k]} y$, and fix $n > k$ sufficiently large that $s_n \mathrel{(\F[k+1][n]
  [n] \setminus \F[k][n][n])} t_n$, where $s_n$ and $t_n$ are the projections of $x$ and $y$ onto 
  $\CantorCantorspace[n][n]$. Then there is a continuous configuration $\gamma_n$ generically 
  compatible with $a_n$ and $\sequence{T_{m,n}}[m, n \in \N]$. Fix $\pair{x_n}{y_n} \in \domain
  {\psi_k^{\gamma_n}}$ with the property that the projections of $x_n$ and $y_n$ onto 
  $\CantorCantorspace[n][n]$ are $s$ and $t$; $\phi^{\gamma_n}(x_n)$, $\phi^{\gamma_n}(y_n)$, 
  $\psi_k^{\gamma_n}(x_n, y_n)$ are extensions of $\phi^{a_n}(s_n)$, $\phi^{a_n}(t_n)$, and 
  $\psi_k^{a_n}(s_n, t_n)$; and $\pair{\pair{\phi^{\gamma_n}(x_n)}{\phi^{\gamma_n}(y_n)}}
  {\psi_k^{\gamma_n}(x_n, y_n)} \in \branches{T_{f(k), f(k+1)}}$. In particular, it follows that $\pair
  {\pair{\phi^{a_n}(s_n)}{\phi^{a_n}(t_n)}}{\psi_k^{a_n}(s_n, t_n)} \in T_{f(k), f(k+1)}$, so $\pair{\pair
  {\phi(x)}{\phi(y)}}{\psi_k(x, y)} \in \branches{T_{f(k), f(k+1)}}$, from which it follows that $\phi(x) 
  \mathrel{(S_{f(k+1)} \setminus R_{f(k)})} \phi(y)$.
\end{theoremproof}

As a corollary, we obtain the following.

\begin{theorem} \label{theorem:closeddichotomy}
  Suppose that $X$ is a \Polish space, $E$ is an analytic equivalence relation on $X$, and 
  $\sequence{R_n}[n \in \N]$ is an increasing sequence of reflexive symmetric $\Fsigma$ binary
  relations on $X$ such that $E \subseteq \union[n \in \N][R_n]$ and \heightcorrection
  {$\relationpower{R_n}{2} \subseteq R_{n+1}$} for all $n \in \N$. Then exactly one of the following 
  holds:
  \begin{enumerate}
    \item The set $X$ is a countable union of \Borel sets on which $\sequence
      {\setcomplement{R_n}}[n \in \N]$ eventually has countable $E$-local \Borel chromatic number.
    \item There exists $f \from \N \to \N$ for which there is a continuous homomorphism $\phi \from
      \CantorCantorspace \to X$ from $\sequence{\F[n], \setcomplement{\F[n]}}[n \in \N]$ to
      $$\sequence{E \intersection R_{f(n)}, \setcomplement{R_{f(n)}}}[n \in \N].$$
  \end{enumerate}
\end{theorem}

\begin{theoremproof}
  In light of Theorem \ref{theorem:Boreldichotomy}, it is sufficient to show that if there is a
  continuous homomorphism $\phi \from \CantorCantorspace \to X$ from $\sequence{\F[n+1] 
  \setminus \F[n]}[n \in \N]$ into $\sequence{R_{n+1} \setminus R_n}[n \in \N]$, then there is a 
  continuous homomorphism from $\sequence{\F[n], \setcomplement{\F[n]}}[n \in \N]$ to 
  $\sequence{E \intersection R_n, \setcomplement{R_n}}[n \in \N]$. Towards this end, define $E' = 
  \preimage{(\phi \times \phi)}{E}$ and $R_n' = \preimage{(\phi \times \phi)}{R_n}$. As Proposition 
  \ref{proposition:containmentmeager} ensures that $R_n'$ is $n$-meager, Proposition \ref
  {proposition:longlonggeneralizedMycielski} yields a continuous homomorphism $\psi \from 
  \CantorCantorspace \to \CantorCantorspace$ from $\sequence{\F[n], \setcomplement{\F[n]}}[n \in 
  \N]$ to $\sequence{\F[n], \setcomplement{R_n'}}[n \in \N]$, in which case the function $\pi = \phi 
  \composition \psi$ is a continuous homomorphism from $\sequence{\F[n], \setcomplement{\F[n]}}
  [n \in \N]$ to $\sequence{E \intersection R_{f(n)}, \setcomplement{R_{f(n)}}}[n \in \N]$.
\end{theoremproof}

\section{Hypersmooth equivalence relations} \label{section:hypersmooth}

In this section, we give a classical proof of Theorem [KL97, Theorem 1]. We 
first note that for witnesses to hypersmoothness, the $\sigma$-ideal appearing in 
Theorems \ref{theorem:Boreldichotomy} and \ref{theorem:closeddichotomy} has a much nicer
characterization.

\begin{proposition} \label{proposition:hypersmooth}
  Suppose that $X$ is a \Polish space, $E$ is a \Borel equivalence relation on $X$, 
  $\sequence{E_n}[n \in \N]$ is an increasing sequence of smooth \Borel equivalence relations on 
  $X$ whose union is $E$, and there are \Borel sets $B_n \subseteq X$, on which 
  $\setcomplement{E_n}$ has countable $E$-local chromatic number, with $X = \union[n \in 
  \N][B_n]$. Then $E$ is essentially hyperfinite.
\end{proposition}

\begin{propositionproof}
  Set $C_n = \union[m < n][B_m]$ and $D_n = B_n \setminus C_n$, and let $F_n$ denote the
  equivalence relation on $X$ given by
  \begin{equation*}
    x \mathrel{F_n} y \Leftrightarrow (x, y \in C_n \mathand x \mathrel{E_n} y) \mathor \exists m \ge n \ (x, y \in 
      D_m \mathand x \mathrel{E_m} y).
  \end{equation*}
  Then $\sequence{F_n}[n \in \N]$ is again an increasing sequence of smooth \Borel equivalence 
  relations whose union is $E$. In addition, $E$ has countable index over $F_0$. Fix \Borel 
  reductions $\phi_m \from X \to \Cantorspace$ of $F_m$ to the equality relation on
  $\Cantorspace$, and observe that the product $\phi \from X \to \CantorCantorspace$, given by 
  $\phi(x)(n) = \phi_n(x)$, is a \Borel reduction of $E$ to $\Eone$. Then $A = \image{\phi}{X}$ is an 
  analytic set on which $\Eone$ is countable, so Theorem \ref{theorem:reflection} yields a \Borel 
  set $B \supseteq A$ on which $\Eone$ is countable. Theorem \ref
  {theorem:countablehypersmooth} then ensures that $\restriction{\Eone}{B}$ is hyperfinite, thus 
  $E$ is essentially hyperfinite.
\end{propositionproof}

As a corollary, we obtain a classical proof of [KL97, Theorem 1].

\begin{theorem}[\Kechris-\Louveau]
  Suppose that $X$ is a \Polish space and $E$ is a hypersmooth \Borel equivalence relation on 
  $X$. Then exactly one of the following holds:
  \begin{enumerate}
    \item The equivalence relation $E$ is essentially hyperfinite.
    \item There is a continuous embedding $\phi \from \CantorCantorspace \to X$ of $\Eone$ into 
      $E$.
  \end{enumerate}
\end{theorem}

\begin{theoremproof}
  Propositions \ref{proposition:meager} and \ref{proposition:essentiallycountable} ensure that the
  two conditions are mutually exclusive.
  
  To see that at least one of them holds, fix an increasing sequence $\sequence{E_n}[n \in \N]$ of 
  smooth \Borel equivalence relations on $X$ whose union is $E$. By Proposition \ref
  {proposition:changeoftopology}, we can assume that each $E_n$ is closed, in which case
  Theorem \ref{theorem:closeddichotomy} and Proposition \ref{proposition:hypersmooth} therefore 
  yield the desired result.
\end{theoremproof}

\section{Treeable equivalence relations} \label{section:treeable}

In this section, we establish our dichotomy theorems for treeable \Borel equivalence relations. 
Given a binary relation $R$ on a set $Y \subseteq X$, we say that a set $Z \subseteq X$ is 
\definedterm{$R$-complete} if $\forall y \in Y \exists z \in Z \ y \mathrel{R} z$.

\begin{proposition} \label{proposition:treeable}
  Suppose that $X$ is a \Polish space, $A \subseteq X$ is analytic, $E$ is a \Borel equivalence
  relation on $X$, $G$ is a \Borel treeing of $E$, and $n$ is a natural number such that for all $x 
  \in A$, there is a countable set $C \subseteq \equivalenceclass{x}{E}$ which is complete with 
  respect to $\restriction{\relationpower{G}{\le n}}{\equivalenceclass{x}{\restriction{E}{A}}}$. Then 
  there is a \heightcorrection{$(\restriction{\relationpower{G}{\le n}}{A})$}-complete \Borel set $B 
  \subseteq X$ on which $E$ is countable.
\end{proposition}

\begin{propositionproof}
  We proceed via induction on $n$. The base case $n = 0$ is trivial, so suppose that we have
  already established the proposition at some $n \in \N$, and for all $x \in A$, there is a countable
  set $C \subseteq \equivalenceclass{x}{E}$ which is complete with respect to 
  $\restriction{\relationpower{G}{\le n + 1}}{\equivalenceclass{x}{\restriction{E}{A}}}$. Let $A'$
  denote the set of $x \in X$ for which there are uncountably many $y \in \verticalsection{G}{x}$ 
  such that for some $m > n$ there is an injective $G$-path $\sequence{z_i}[i \le m]$ with $x = 
  z_0$, $y = z_1$, and $z_m \in A$. As Theorem \ref{theorem:uncountablequantifier} ensures 
  that the property of being countable is $\Piclass[1][1]$-on-$\Sigmaclass[1][1]$, the set $A'$ is 
  analytic. Moreover, the acyclicity of $G$ ensures that if $x \in A$ and $C \subseteq 
  \equivalenceclass{x}{E}$ is a countable set which is complete with respect to $\restriction
  {\relationpower{G}{\le n + 1}}{\equivalenceclass{x}{\restriction{E}{A}}}$, then $A' \intersection 
  \equivalenceclass{x}{E} \subseteq C$. In particular, it follows that $E$ is countable on $A'$. As 
  this latter property is again $\Piclass[1][1]$-on-$\Sigmaclass[1][1]$, Theorem \ref
  {theorem:reflection} yields a \Borel set $B' \supseteq A'$ on which $E$ is countable. As Theorems 
  \ref{theorem:setofunicity} and \ref{theorem:bianalytic} ensure that $\relationpower{G}{\le n+1}$ is 
  \Borel, Theorem \ref{theorem:countabletoone} implies that the set $B''$ of points $\relationpower
  {G}{\le n +1}$-related to points in $B'$ is \Borel.
  
  Define $A'' = A \setminus B''$, and observe that if $x \in A''$ and $C \subseteq \equivalenceclass
  {x}{E}$ is a countable set which is complete with respect to $\restriction
  {\relationpower{G}{\le n + 1}}{\equivalenceclass{x}{\restriction{E}{A}}}$, then there exists $y \in C 
  \setminus B'$ such that $x$ is $\relationpower{G}{\le n}$-related to either $y$ or one of its 
  countably many neighbors $z$ for which $\pair{y}{z}$ extends to a $G$-path $\sequence{w_i}[i 
  \le n + 1]$ from $y$ to $A$. In particular, the induction hypothesis yields a $(\restriction
  {\relationpower{G}{\le n}}{A''})$-complete \Borel set $B''' \subseteq X$ on which $E$ is countable, 
  in which case the set $B = B' \union B'''$ is as desired.
\end{propositionproof}

As corollaries, we obtain the following dichotomy theorems.

\begin{theorem} \label{theorem:homomorphismtreeabledichotomy}
  Suppose that $X$ is a \Polish space, $E$ is a treeable \Borel equivalence relation on $X$, and 
  $G$ is a \Borel treeing of $E$. Then exactly one of the following holds:
  \begin{enumerate}
    \item There is an $E$-complete \Borel set on which $E$ is countable.
    \item There exists a function $f \from \N \to \N$ for which there is a continuous homomorphism
      $\phi \from \CantorCantorspace \to X$ from $\sequence{\F[n+1] \setminus \F[n]}[n \in \N]$ to 
      $\sequence{\relationpower{G}{\le f(n+1)} \setminus \relationpower{G}{\le f(n)}}[n \in \N]$.
  \end{enumerate}
\end{theorem}

\begin{theoremproof}
  As condition (1) ensures that $X$ is of the form $\union[n \in \N][B_n]$, where each $B_n 
  \subseteq X$ is a \Borel set on which $\setcomplement{\relationpower{G}{\le n}}$ has countable 
  $E$-local \Borel chromatic number, Proposition \ref{proposition:meager} ensures that the two 
  conditions are mutually exclusive. Theorem \ref{theorem:Boreldichotomy} and Proposition \ref
  {proposition:treeable} imply that at least one of them holds.  
\end{theoremproof}

\begin{theorem} \label{theorem:originaltreeabledichotomy}
  Suppose that $X$ is a \Polish space and $E$ is a \Borel equivalence relation on $X$ which is
  subtreeable-with-$\Fsigma$-iterates. Then for every analytic set $A \subseteq X$, exactly one of 
  the following holds:
  \begin{enumerate}
    \item There is an $(\restriction{E}{A})$-complete \Borel set $B \subseteq X$ on which $E$ is 
      countable.
    \item There is a continuous embedding $\phi \from \CantorCantorspace \to X$ of $\Eone$ into
      $\restriction{E}{A}$.
  \end{enumerate}
\end{theorem}

\begin{theoremproof}
  Proposition \ref{proposition:meager} ensures that the two conditions are mutually exclusive, and 
  Theorem \ref{theorem:closeddichotomy} and Proposition \ref{proposition:treeable} imply that at 
  least one of them holds.
\end{theoremproof}

We say that \emph{embeddability of $\Eone$ is determined below $E$ by $\calE$} if for every 
analytic set $A \subseteq X$, either $\restriction{E}{A} \in \calE$ or there is a continuous 
embedding of $\Eone$ into $E$. Theorem \ref{theorem:originaltreeabledichotomy} implies \Borel 
equivalence relations which are subtreeable-with-$\Fsigma$-iterates have this property, where
$\calE$ is the class of essentially countable \Borel equivalence relations on \Polish spaces. The
following fact implies that this holds under the weaker assumption of being essentially 
subtreeable-with-$\Fsigma$-iterates.

\begin{proposition} \label{proposition:treeabledichotomy}
  Suppose that $\calE$ is a class of \Borel equivalence relations on \Polish spaces. Then the 
  class of \Borel equivalence relations below which embeddability of $\Eone$ is determined by 
  essentially $\calE$ is closed under \Borel reducibility.
\end{proposition}

\begin{propositionproof}
  Suppose that $X$ and $Y$ are \Polish spaces, $E$ and $F$ are \Borel equivalence relations on
  $X$ and $Y$, $\pi \from X \to Y$ is a \Borel reduction of $E$ to $F$, and embeddability of
  $\Eone$ is determined below $F$ by essentially $\calE$. Given an analytic set $A \subseteq
  X$, either there is a \Borel reduction $\psi$ of $\restriction{F}{\image{\pi}{A}}$, and therefore of 
  $\restriction{E}{A}$, to a \Borel equivalence relation in $\calE$, or there is a continuous 
  embedding of $\Eone$ into $\restriction{F}{\image{\pi}{A}}$, in which case Proposition \ref
  {proposition:Eoneclosure} yields a continuous embedding of $\Eone$ into $\restriction{E}{A}$. 
\end{propositionproof}

\section{Borel functions} \label{section:Kakutani}

In this section, we establish a natural strengthening of Theorem \ref
{theorem:originaltreeabledichotomy} for graphs induced by functions. Although not strictly 
necessary to achieve this goal, we will first establish several preliminary results so as to further 
clarify the nature of essential countability in this context.

\begin{proposition} \label{proposition:essentiallycountablereflection}
  Suppose that $X$ and $Y$ are \Polish spaces, $E$ is a \Borel equivalence relation on $X$, $F$
  is a countable equivalence relation on a subset of $Y$, and $\pi \from X \to Y$ is a \Borel
  reduction of $E$ to $F$. Then there is a countable \Borel equivalence relation $F'$ on $Y$ such
  that $\pi$ is also a reduction of $E$ to $F'$.
\end{proposition}

\begin{propositionproof}
  Set $R = \image{(\pi \times \pi)}{E}$. The fact that $\pi$ is a homomorphism from $E$ to $F$
  ensures that $R \subseteq F$. As $F$ is countable and $\pi$ is a cohomomorphism from $E$ to
  $F$, it follows that $R$ is subset of $Y \times Y$, with countable horizontal and vertical 
  sections, for which $\pi$ is a cohomomorphism from $E$ to the smallest equivalence relation on
  $Y$ containing $R$. As $R$ is analytic and this latter property is $\Piclass[1]
  [1]$-on-$\Sigmaclass[1][1]$, Theorem \ref{theorem:reflection} yields a \Borel set $R' \supseteq 
  R$, with countable horizontal and vertical sections, for which $\pi$ is a cohomomorphism from 
  $E$ to the smallest equivalence relation on $Y$ containing $R'$. Let $F'$ denote the latter 
  equivalence relation. As $R \subseteq F'$, it follows that $\pi$ is also a homomorphism from $E$ 
  to $F'$, and therefore $\pi$ is a reduction of $E$ to $F'$. As the horizontal and vertical sections of 
  $R'$ are countable, it follows that $F'$ is countable, so Theorem \ref{theorem:countabletoone}
  ensures that $F'$ is \Borel.
\end{propositionproof}

\begin{proposition} \label{proposition:countableindex}
  Suppose that $X$ is a \Polish space, $E$ and $F$ are \Borel equivalence relations on $X$, 
  and $E \intersection F$ has countable index in $E$ and $F$. Then $E \text{ is essentially 
  countable} \iff F \text{ is essentially countable}$.
\end{proposition}

\begin{propositionproof}
  It is sufficient to handle the special case that $E \subseteq F$.
  
  To see $(\implies)$, suppose that $X'$ is a \Polish space, $E'$ is a countable equivalence 
  relation on $X'$, and $\pi \from X \to X'$ is a \Borel reduction of $E$ to $E'$. Then $\pi$ is a 
  reduction of $F$ to the countable equivalence relation $\image{(\pi \times \pi)}{F}$ on $\image
  {\pi}{X}$, so $F$ is essentially countable by Proposition \ref
  {proposition:essentiallycountablereflection}.
  
  To see $(\impliedby)$, suppose that $X'$ is a \Polish space, $F'$ is a countable equivalence 
  relation on $X'$, and $\phi \from X \to X'$ is a \Borel reduction of $F$ to $F'$. Let $D'$ denote the 
  equality relation on $X'$, and observe that the relation $D = \preimage{(\phi \times \phi)}{D'}$ has 
  countable index in $F$. By $(\implies)$, it is enough to show that $D \intersection E$ is smooth, 
  thus essentially countable. Suppose, towards a contradiction, that this is not the case. Then 
  Theorem \ref{theorem:E0dichotomy} yields a continuous embedding $\psi \from 
  \Cantorspace \to X$ of $\Ezero$ into $D \intersection E$, in which case $\phi \composition \psi$ is 
  a countable-to-one \Borel homomorphism from $\Ezero$ to $D'$, contradicting Proposition
  \ref{proposition:meagertoonehomomorphism}.
\end{propositionproof}

With these preliminaries out of the way, we now turn our attention to functions $T \from X \to
X$. Let \definedsymbol{\Etail[T]} denote the equivalence relation on $X$ given by $x \mathrel
{\Etail[T]} y \iff \exists m, n \in \N \ T^m(x) = T^n(y)$.

The \definedterm{eventually periodic part} of $T$ is the set of $x \in X$ for which there are
natural numbers $m < n$ with $T^m(x) = T^n(x)$, and $T$ is \definedterm{aperiodic} if its
eventually periodic part is empty. The following observation will allow us to focus our attention
on aperiodic functions.

\begin{proposition} \label{proposition:periodicpart}
  Suppose that $X$ is a \Polish space and $T \from X \to X$ is \Borel.
  Then there is a \Borel transversal of the restriction of $\Etail[T]$ to the eventually periodic part
  of $T$.
\end{proposition}

\begin{propositionproof}
  The \definedterm{periodic part} of $T$ is the set of $x \in X$ for which there is a positive natural
  number $n$ with $x = T^n(x)$. As the periodic part of $T$ intersects every equivalence class of
  $\Etail[T]$ in a finite set, the desired result follows from the fact that every finite \Borel equivalence
  relation on a \Polish space has a \Borel transversal, which itself is a consequence of 
  Theorem \ref{theorem:LusinNovikov}.
\end{propositionproof}

The following observation will allow us to apply our earlier results.

\begin{proposition} \label{proposition:functiontreeability}
  Suppose that $X$ is a \Polish space and $T \from X \to X$ is \Borel. Then
  $\Etail[T]$ is treeable.
\end{proposition}

\begin{propositionproof}
  As \Borel equivalence relations with \Borel transversals are trivially treeable, Proposition \ref
  {proposition:periodicpart} allows us to assume that $T$ is aperiodic. Then the
  graph \definedsymbol{\functiongraph{T}} on $X$ given by $x \mathrel{\functiongraph{T}} y \iff 
  (T(x) = y \mathor T(y) = x)$ is a \Borel treeing of $E$.
\end{propositionproof}

This yields another characterization of essential countability of $$\Etail[T].$$

\begin{proposition} \label{proposition:essentiallycountableequivalents}
  Suppose that $X$ is a \Polish space, $T \from X \to X$ is \Borel, and
  $\Etail[T]$ is essentially countable. Then there is an $\Etail[T]$-complete \Borel set on which
  $\Etail[T]$ is countable.
\end{proposition}

\begin{propositionproof}
  By Proposition \ref{proposition:functiontreeability}, the equivalence relation $\Etail[T]$ is treeable.
  The desired result is therefore a consequence of Theorem \ref
  {theorem:treeableessentiallycountable}. Although this latter result has a classical proof (see
  [Mil12]), we will give a simpler argument using the structure of $T$.
  
  By Proposition \ref{proposition:periodicpart}, we can assume that $T$ is aperiodic. By Proposition
  \ref{proposition:changeoftopology}, we can assume that $X$ carries a \Polish topology with 
  respect to which $T$ is continuous. Then the iterates of $\functiongraph{T}$ are closed. 
  Theorem \ref{theorem:originaltreeabledichotomy} therefore yields the desired $\Etail
  [T]$-complete \Borel set on which $\Etail[T]$ is countable.
\end{propositionproof}

Define \definedsymbol{\Ezero[T]} on $X$ by $x \mathrel{\Ezero[T]} y \iff \exists n \in \N \ T^n(x) = 
T^n(y)$. Note that $\Ezero(T)$ is a countable index subequivalence relation of $\Etail(T)$.

\begin{proposition}
  Suppose that $X$ is a \Polish space and $T \from X \to X$ is \Borel. Then
  $\Ezero[T]$ is essentially countable if and only if $\Etail[T]$ is essentially countable.
\end{proposition}

\begin{propositionproof}
  This is a direct consequence of Proposition \ref{proposition:countableindex}.
\end{propositionproof}

Together with Proposition \ref{proposition:essentiallycountableequivalents}, the following 
fact ensures that $\Etail[T]$ is essentially countable if and only if $T$ is essentially 
countable-to-one.

\begin{proposition} \label{proposition:essentiallycountabletoone}
  Suppose that $X$ is a \Polish space, $T \from X \to X$ is \Borel, and $B \subseteq X$ is a \Borel 
  set on which $\Etail[T]$ is countable. Then there is a $T$-stable \Borel set $A \supseteq B$ on 
  which $\Etail[T]$ is countable.
\end{proposition}

\begin{propositionproof}
  Set $A = \union[n \in \N][\image{T^n}{B}]$. Then $A$ is $T$-stable, and Theorem \ref
  {theorem:countabletoone} ensures that it is \Borel.
\end{propositionproof}

Define \definedsymbol{\F[n](T)} on $X$ by $x \mathrel{\F[n](T)} y \iff T^n(x) = T^n(y)$.

\begin{proposition} \label{proposition:kakutanione}
  Suppose that $X$ is a \Polish space, $T \from X \to X$ is \Borel, and there is a sequence 
  $\sequence{B_n}[n \in \N]$ of \Borel sets for which $X = \union[n \in \N][B_n]$ and $\sequence
  {\setcomplement{\F[k](T)}}[k \in \N]$ eventually has countable $\Ezero(T)$-local chromatic 
  number for all $n \in \N$. Then $T$ is essentially countable-to-one.
\end{proposition}

\begin{propositionproof}
  Fix natural numbers $k_n \in \N$ such that $\setcomplement{\F[k_n](T)}$ has countable
  $\Ezero[T]$-local chromatic number on $B_n$ for all $n  \in \N$. Then $\Ezero[T]$ is countable 
  on the analytic set $A = \union[n \in \N][\image{T^{k_n}}{B_n}]$. As Theorem \ref
  {theorem:uncountablequantifier} ensures that the property of being countable is $\Piclass[1]
  [1]$-on-$\Sigmaclass[1][1]$, Theorem \ref{theorem:reflection} yields a \Borel set $B \supseteq A$ 
  on which $\Ezero[T]$ is countable. As $\Etail[T]$ must also be countable on this set, Proposition
  \ref{proposition:essentiallycountabletoone} ensures that $T$ is essentially countable-to-one.
\end{propositionproof}

Define \definedsymbol{\R[n][T]} on $X$ by $x \mathrel{\R[n][T]} y \iff \exists i, j \le n \ T^i(x) = T^j(y)$.

\begin{proposition} \label{proposition:kakutanitwo}
  Suppose that $X$ is a \Polish space, $T \from X \to X$ is aperiodic, $f \from \N \to \N$, and $\phi 
  \from \CantorCantorspace \to X$ is a homomorphism from $\sequence{\F[n], \setcomplement{\F
  [n]}}[n \in \N]$ into $\sequence{\F[f(n)](T), \setcomplement{\R[f(n)][T]}}[n \in \N]$. Then the function
  $\pi \from \N \times \CantorCantorspace \to X$ given by $\pi(n, s^n(x)) = T^{f(n)} \composition 
  \phi(x)$ defines a \Kakutani embedding of $S \times s$ into $T$.
\end{proposition}

\begin{propositionproof}
  To see that $\pi$ is well-defined, note that if $s^n(x) = s^n(y)$, then $x \mathrel{\F[n]} y$, so $\phi
  (x) \mathrel{\F[f(n)](T)} \phi(y)$, thus $T^{f(n)} \composition \phi(x) = T^{f(n)} \composition \phi(y)$.
  Note also that the set $B = \image{\pi}{\N \times \CantorCantorspace}$ is trivially \definedterm
  {$T$-recurrent}, in the sense that $B \subseteq \union[n > 0][\image{T^{-n}}{B}]$.
  
  To see that $\pi$ is injective, suppose that $m, n \in \N$ and $x, y \in \CantorCantorspace$ are
  such that $\pi(m, x) = \pi(n, y)$. By reversing the roles of $x$ and $y$ if necessary, we can
  assume that $m \le n$. Fix $x', y' \in \CantorCantorspace$ such that $x = s^m(x')$ and $y = s^n
  (y')$, and observe that $T^{f(m)} \composition \phi(x') = \pi(m, x) = \pi(n, y) = T^{f(n)} \composition
  \phi(y')$, so the fact that $f(m) \le f(n)$ ensures that $\phi(x') \mathrel{\R[f(n)][T]} \phi(y')$. As $\phi$ 
  is a homomorphism from $\setcomplement{\F[n]}$ to $\setcomplement{\R[f(n)][T]}$, it follows that 
  $x' \mathrel{\F[n]} y'$. As
  $\phi$ is also a homomorphism from $\F[n]$ to $\F[f(n)](T)$, it follows that $T^{f(n)} \composition
  \phi(x') = T^{f(n)} \composition \phi(y')$. Then $T^{f(m)} \composition \phi(x') = T^{f(n)} 
  \composition \phi(x')$, so the injectivity of $f$ and the aperiodicity of $T$ ensure that 
  $m = n$, thus $x = s^n(x') = s^n(y') = y$.
  
  Suppose now that $n \in \N$ and $x \in \CantorCantorspace$, and fix $x' \in \CantorCantorspace$
  for which $x = s^n(x')$. As $\phi$ is a homomorphism from $\pair{\Ezero(s)}{\setcomplement
  {\Ezero(s)}}$ to $\pair{\Ezero(T)}{\setcomplement{\Etail(T)}}$, it follows that $\image{\phi}
  {\equivalenceclass{x'}{\Ezero(s)}} = \image{\phi}{\CantorCantorspace} \intersection 
  \equivalenceclass{\phi(x')}{\Ezero(T)} = \image{\phi}{\CantorCantorspace} \intersection 
  \equivalenceclass{\phi(x')}{\Etail(T)}$, thus
  \begin{align*}
    T_{\image{\phi}{\N \times \CantorCantorspace}} \composition \pi(n, x) 
      & = T_{\image{\phi}{\N \times \CantorCantorspace}} \composition \pi(n, s^n(x')) \\
      & = T_{\image{\phi}{\N \times \CantorCantorspace}} \composition T^{f(n)} \composition \phi(x') \\
      & = T^{f(n+1)} \composition \phi(x') \\
      & = \pi(n+1, s^{n+1}(x')) \\
      & = \pi((S \times s)(n, x)), \\
  \end{align*}
  thus $\pi$ is a \Kakutani embedding of $S \times s$ into $T$.
\end{propositionproof}

We are now ready to establish our final result.

\begin{theorem}
  Suppose that $X$ is a \Polish space and $T \from X \to X$ is \Borel. Then exactly one of 
  the following holds:
  \begin{enumerate}
    \item The function $T$ is essentially countable-to-one.
    \item There is a continuous \Kakutani embedding $\phi \from \CantorCantorspace \to X$ of $S 
      \times s$ into $T$.
  \end{enumerate}
\end{theorem}

\begin{theoremproof}
  To see that the two conditions are mutually exclusive, suppose that $B$ is a $T$-complete, 
  $T$-stable \Borel set on which $T$ is countable-to-one, and $\pi \from \N \times 
  \CantorCantorspace \to X$ is a \Borel \Kakutani embedding of $S \times s$ into $T$. Then 
  $\preimage{\pi}{B}$ is an $(S \times s)$-complete, $(S \times s)$-stable \Borel set on which $S 
  \times s$ is countable-to-one, so $\image{\projection{\CantorCantorspace}}{\preimage{\pi}{B}}$ 
  is an $\Eone$-complete \Borel set on which $\Eone$ is countable, contradicting Proposition \ref
  {proposition:meager}.

  It remains to check that at least one of the two conditions holds. By Proposition \ref
  {proposition:periodicpart}, we can assume that $T$ is aperiodic. By Proposition \ref
  {proposition:changeoftopology}, we can assume that $T$ is continuous, in which case each of
  the relations $\R[n][T]$ is closed. Theorem \ref{theorem:closeddichotomy} and Proposition
  \ref{proposition:kakutanione} ensure that if $T$ is not essentially countable-to-one, then there is
  a function $f \from \N \to \N$ for which there is a continuous homomorphism from $\sequence{\F
  [n], \setcomplement{\F[n]}}[n \in \N]$ to $\sequence{\Ezero[T] \intersection \R[f(n)][T], 
  \setcomplement{\R[f(n)][T]}}[n \in \N]$. As the aperiodicity of $T$ implies that $\F[n](T) = \Ezero[T] 
  \intersection \R[n][T]$ for all $n \in \N$, Proposition \ref{proposition:kakutanitwo} yields a 
  continuous \Kakutani embedding of $S \times s$ into $T$.
\end{theoremproof}

\begin{acknowledgements}
  We would like to thank Clinton Conley for allowing us to include Proposition \ref
  {proposition:perfectpair} and its proof, as well as Clinton Conley, Alexander Kechris, and the
  anonymous referee for their comments on earlier versions of the paper.
\end{acknowledgements}

\section{References}

\noindent [DJK94] R. Dougherty, S. Jackson and A. S. Kechris,\it\ The structure of hyperfinite Borel equivalence relations,\rm\ Trans. Amer. Math. Soc. 341, 1 (1994), 193-225 

\noindent [Hjo08] G. Hjorth,\it\ Selection theorems and treeability,\rm\ Proc. Amer. Math. Soc. 136, 10 (2008), 3647-3653

\noindent [HKL90] L. A. Harrington, A. S. Kechris and A. Louveau,\it\ A Glimm-Effros dichotomy for Borel equivalence relations,\rm\ J. Amer. Math. Soc. 3, 4 (1990), 903-928

\noindent [Kec95] A. S. Kechris,\it\ Classical descriptive set theory,\rm\ Graduate Texts in Mathematics, vol. 156, Springer-Verlag, New York (1995)

\noindent [KL97] A. S. Kechris and A. Louveau,\it\ The classification of hypersmooth Borel equivalence relations,\rm\ J. Amer. Math. Soc. 10, 1 (1997), 215-242

\noindent [KST99] A. S. Kechris, S. Solecki and S. Todorcevic,\it\ Borel chromatic numbers,\rm\ Adv. Math. 141, 1 (1999), 1-44

\noindent [Kur68] K. Kuratowski,\it\ Topology. Vol. II,\rm\ New edition, revised and augmented. Translated from the French by A. Kirkor, Academic Press, New York (1968)

\noindent [Mil12] B. D. Miller,\it\ The graph-theoretic approach to descriptive set theory,\rm\ Bull. Symbolic Logic 18, 4 (2012), 554-575

\noindent [Sil80] J. H. Silver,\it\ Counting the number of equivalence classes of Borel and coanalytic equivalence relations,\rm\ Ann. Math. Logic 18, 1 (1980), 1-28

\end{document}